\documentclass[%
 amsmath,amssymb,
preprint,%
10pt,
a4paper,
twocolumn,
]{revtex4-1}

\usepackage{graphicx}% Include figure files
\usepackage{dcolumn}% Align table columns on decimal point
\usepackage{bm}% bold math
\usepackage[utf8]{inputenc}
\usepackage[T1]{fontenc}
\usepackage{mathptmx}
\usepackage{etoolbox}
\usepackage[dvipsnames]{xcolor}

\makeatletter
\def\@email#1#2{%
 \endgroup
 \patchcmd{\titleblock@produce}
  {\frontmatter@RRAPformat}
  {\frontmatter@RRAPformat{\produce@RRAP{*#1\href{mailto:#2}{#2}}}\frontmatter@RRAPformat}
  {}{}
}%
\makeatother

\graphicspath{{Figures/}}

\begin{document}

%\preprint{KEA23}

\title[Transitional Cluster States in Delay-Coupled Oscillators]{
Transitional cluster dynamics in a model for delay-coupled chemical oscillators}
% Force line breaks with \\
\author{Andrew Keane}
 \email{andrew.keane@ucc.ie}
\affiliation{ 
School of Mathematical Sciences, University College Cork, Cork, T12 XF62, Ireland%\\This line break forced with \textbackslash\textbackslash
}%
\affiliation{ 
Environmental Research Institute, University College Cork, Cork, T23 XE10, Ireland%\\This line break forced with \textbackslash\textbackslash
}%
\author{Alannah Neff}%
% \email{118369773@umail.ucc.ie}
\affiliation{ 
School of Mathematical Sciences, University College Cork, Cork, T12 XF62, Ireland%\\This line break forced with \textbackslash\textbackslash
}%

\author{Karen Blaha}
% \email{kb4hk@virginia.edu}
\affiliation{%
Sandia National Labs, 1515 Eubank Blvd SE1515 Eubank Blvd SE, Albuquerque, NM 87123, USA %\\This line break forced% with \\
}%

\author{Andreas Amann}%
% \email{a.amann@ucc.ie}
\affiliation{ 
School of Mathematical Sciences, University College Cork, Cork, T12 XF62, Ireland%\\This line break forced with \textbackslash\textbackslash
}%

\author{Philipp H\"ovel}
% \email{philipp.hoevel@gmail.com}
\affiliation{%
Department of Electrical and Information Engineering, Christian-Albrechts-Universit\"at zu Kiel, Kaiserstr. 2, 24143 Kiel, Germany%\\This line break forced% with \\
}%

\collaboration{in memory of John L. Hudson}

\date{May 2023}% It is always \today, today,

\begin{abstract}

Cluster synchronization is a fundamental phenomenon in systems of coupled oscillators. Here, we investigate clustering patterns that emerge in a unidirectional ring of four delay-coupled electrochemical oscillators. A voltage parameter in the experimental set-up controls the onset of oscillations via a Hopf bifurcation. For a smaller voltage, the oscillators exhibit simple, so-called primary, clustering patterns, where all phase differences between each set of coupled oscillators are identical. However, upon increasing the voltage, additional secondary states, where phase differences differ, are detected. Previous work on this system saw the development of a mathematical model that explains how the existence, stability, and common frequency of the experimentally observed cluster states can be accurately controlled by the delay time of the coupling.

In this study, we revisit the mathematical model of the electrochemical oscillators to address open questions by means of bifurcation analysis. Our analysis reveals how the stable cluster states, corresponding to experimental observations, lose their stability via an assortment of bifurcation types. The analysis further reveals a complex interconnectedness between branches of different cluster types; in particular, we find that each secondary state provides a continuous transition between certain primary states.  These connections are explained by studying the phase space and parameter symmetries of the respective states. Furthermore, we show that it is only for a larger value of the voltage parameter that the branches of secondary states develop intervals of stability. Otherwise, for a smaller voltage, all the branches of secondary states are completely unstable and therefore hidden to experimentalists.

\end{abstract}

\maketitle

\section{\label{sec:Intro}Introduction}

Nonlinear systems of coupled oscillators form the basis of multiple areas of interdisciplinary research, from the dynamics of coupled lasers to modelling neuronal dynamics. Of particular importance and wide-spread interest is the phenomenon of synchronization and clustering \cite{juang2014, lodi2020, protachevicz2021, soriano2013, han2019}. For example, in neuronal dynamics different patterns of synchronization are related to normal cognitive and pathological functions of the brain \cite{schnitzler2005}. In power grid dynamics, analyzing synchronization and possible cluster states can provide insights into the stability of the grid \cite{motter2013, pecora14}.

Clustering behavior is often associated with underlying symmetry properties of a system \cite{macarthur08, skardal19}.
For such systems, group theory and equivariant dynamical system theory can provide a useful tool for exploring and understanding the dynamics \cite{chossat00, golubitsky12}. 
For example, the symmetry properties of a network of oscillators can be used to facilitate model reductions \cite{pietras19}.
Knowledge of the symmetries can also be used to discover possible cluster synchronization patterns \cite{pecora14}.
Furthermore, the authors of Ref.~[\onlinecite{nicosia13}] show that nodes in a complex network will form clusters according to their symmetry properties within the network. It is, therefore, suggested that structural connectivity of the brain could play a role in neural synchronization across distant locations. 

There have also been many studies on the interplay between dynamics and symmetry in the context of networks with \emph{ring} topologies, which will be the focus of this paper.
For example, in Refs.~[\onlinecite{schneider13,schneider16}] a type of equivariant delayed feedback control is used to target the stabilization of periodic orbits with a specified spatio-temporal pattern in rings of oscillators. In Refs.~[\onlinecite{collins93a,collins93b}] rings of oscillators are used to study the gaits (i.e. walking patterns) of animals. As a result, it is suggested that transitions between different gaits can be modelled as symmetry-breaking bifurcations: for example, a horse will walk, trot, then gallop, as a result of successive bifurcations.

In more recent years, greater attention has been paid to the fact that many processes are not instantaneous, but instead possess an inherent delay. Such delays can potentially have a crucial influence on the overall dynamics of the system \cite{atay2010, zakharova13, erneux2017, otto2019}. In certain cases, it can turn trivial dynamics into complex dynamics \cite{calleja2017}.

In order to better understand the effects of delay and symmetry on synchronization and clustering, both experimental research and mathematical modelling have an important role to play --- one can inform the other of where interesting dynamical phenomena may occur and provides clues to possible mechanisms behind the phenomena. 
For example, in the context of mutually coupled lasers, mathematical modelling has explained the various dynamical regimes through the existence of mutually synchronized symmetric and symmetry-broken states which are connected through pitchfork bifurcations \cite{yanchuk2004dynamics,erzgraber2006compound,clerkin2014multistabilities}.
Clustering has also been studied in networks of physically dissimilar mechanical and electrical oscillators \cite{della2020symmetries, blaha2016symmetry}, where the interplay between modelling and experiment focus experiments on interesting parametric regions and highlights which dynamics may be most realizable experimentally. 

The work presented here builds upon a previous study on the clustering patterns observed in an experimental set-up of four electrochemical oscillators coupled in a unidirectional ring \cite{BLA13}. It was shown that a delay in the coupling between each set of neighboring oscillators could control which clustering pattern would emerge. Another important parameter of the experimental set-up (see Appendix~\ref{app:experiment} for a brief review) is the choice of applied voltage. It is only above a critical voltage that the oscillators actually oscillate via a Hopf bifurcation. Above this critical value, two voltage regimes were studied. One lower voltage regime, which resulted in smooth sinusoidal oscillations. This was called the \emph{smooth} regime. In the second regime, a higher voltage resulted in relaxation oscillations, and was called the \emph{relaxation} regime.

\begin{figure}[t]
\includegraphics[width=\linewidth]{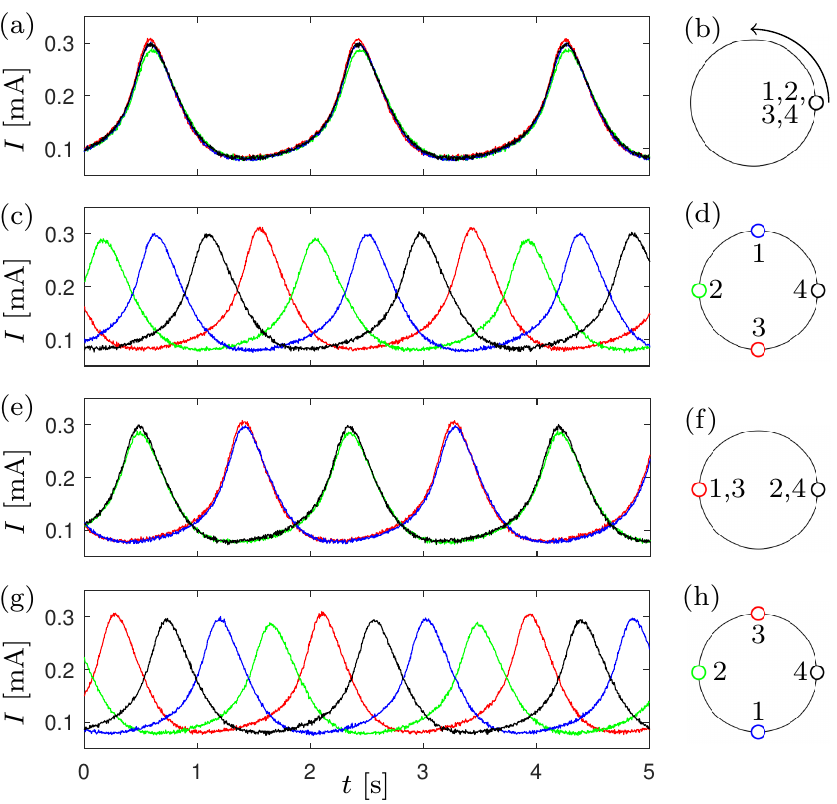}
\caption{Experimental time series and schematic diagrams of primary cluster states in the smooth regime with $V_0 = 1.105$V and $K=0.15$: (a)--(b) in-phase with $\tau= 0.95\times\left(\frac{2\pi}{\omega}\right)$, (c)--(d) splay with $\tau= 1.25\times\left(\frac{2\pi}{\omega}\right)$, (e)--(f) 2-cluster  with $\tau= 0.50\times\left(\frac{2\pi}{\omega}\right)$, and (g)--(h) reverse splay states  with $\tau= 0.70\times\left(\frac{2\pi}{\omega}\right)$. }
\label{fig:exp_primary}
\end{figure}

\begin{figure}[th]
\includegraphics[width=\linewidth]{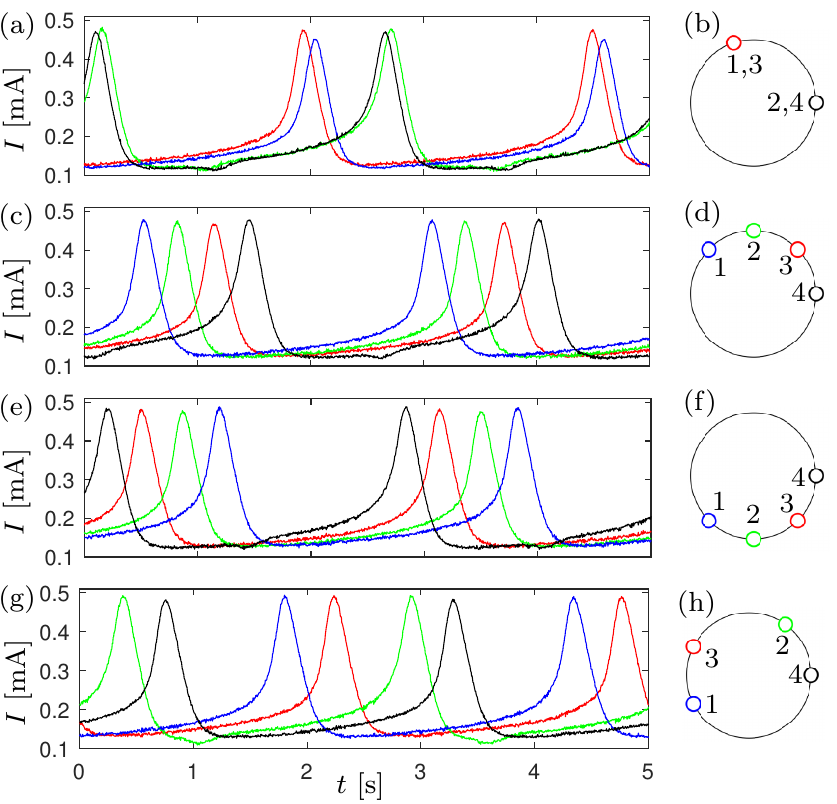}
\caption{Experimental time series and schematic diagrams of secondary cluster states in the relaxation regime with $V_0 = 1.2$V and $K=0.10$: (a)--(b) compressed 2-cluster  with $\tau= 0.65\times\left(\frac{2\pi}{\omega}\right)$, (c)--(d) compressed reverse splay with $\tau= 0.81\times\left(\frac{2\pi}{\omega}\right)$, (e)--(f) compressed splay with $\tau= 1.06\times\left(\frac{2\pi}{\omega}\right)$, and (g)--(h) open 2-cluster states with $\tau= 0.51\times\left(\frac{2\pi}{\omega}\right)$. }
\label{fig:exp_secondary}
\end{figure}

The clustering patterns that were observed depended on the regime. In the smooth regime, only clustering patterns with equal phase differences between neighboring oscillators (neighboring in the sense of ring network topology) were found. These are called the \emph{primary} clustering patterns. 
A clear visualization of the clustering patterns in the smooth and relaxation regimes is crucial for appreciating the experimentally-motivated modelling, as well as reconnecting the theoretical analysis of the model back to the experimental observations. Therefore, we show examples of the different primary clustering patterns in Fig.~\ref{fig:exp_primary}, as both time series and schematic diagrams.
Panels~(a) and (b) show an example of an \emph{in-phase} state, where all oscillators possess equal phase. Panels~(c) and (d) depict an example of the \emph{splay} state, where all neighboring oscillators have a common phase difference of $\varphi_{j+1}-\varphi_j = \pi/2$ with their neighbors. A \emph{2-cluster} state is shown in panels~(e) and (f), where neighboring oscillators have a phase difference of $\varphi_{j+1}-\varphi_j = \pi$. Similarly, panels~(g) and (h) show a \emph{reverse splay} state, with phase differences $\varphi_{j+1}-\varphi_j = 3\pi/2$. 
We calculate phase via peak-to-peak linear interpolation\cite{rusin2010synchronization}.

In the relaxation regime, the primary clustering patterns were still observed, however additional \emph{secondary} clustering patterns were observed with phase differences that were not all equal.
The first example shown in Fig.~\ref{fig:exp_secondary}(a) and (b) is a \emph{compressed 2-cluster} state. It resembles a 2-cluster state, only that the phase difference between the two clusters is less than $\pi$. Panels~(c) and (d) show an example of a \emph{compressed reverse splay} state. While three of the phase differences between neighboring oscillators appear to remain equal, one of the phase differences is less than $3\pi/2$. Similarly, a \emph{compressed splay} state, like the one shown in panels~(e) and (f), appears to maintain three equal phase differences between neighboring oscillators, while one phase difference becomes greater than $\pi/2$. Finally, \emph{open 2-cluster} states are often observed in the experiments, as depicted in panels~(g) and (h). This state resembles a 2-cluster state, where the oscillators within each cluster appear to have drifted apart.

A well-calibrated model of Stuart-Landau oscillators was found to match the transitions between different clustering patterns, as well as their common frequency, in the smooth regime with very good accuracy. 
Furthermore, an extension to the classic Stuart-Landau oscillators with experimentally-derived coupling functions was found to agree remarkably well with the experimental results in the relaxation regime. The results presented in Ref.~[\onlinecite{BLA13}] effectively took the form of bifurcation diagrams in terms of the coupling delay parameter, with the branches of various clustering patterns in the experiment being accurately matched by stable branches of corresponding solutions to the model. However, because the focus of that study was on the accurate reproduction of experimental results, the analysis of the model was limited to only the stable behavior. While successful, this approach left important open questions. It was still unclear how these states lose their stability. Perhaps more intriguing is the question: Where do the secondary states in the relaxation regime come from, when there is no evidence of their existence in the smooth regime?

In this paper, we conduct a systematic bifurcation analysis of the model for the unidirectional ring of delay-coupled electrochemical oscillators that reveals the roles of the unstable solutions and various bifurcation types in the organization of clustering patterns with respect to the coupling delay.
More specifically, we consider each type of secondary cluster state observed in the experiments and identify how they relate to other cluster states. A useful tool in this context is the investigation of the symmetry groups of the various cluster states. We find that each secondary cluster state can be interpreted as a transitional cluster state, en route between certain primary states via symmetry-breaking bifurcations. Primary states are associated with symmetry groups of order 4, while secondary states have symmetry groups of order 2 or 1. Finally, we find that although the general interconnected structure of clustering patterns is preserved in the smooth regime, all secondary cluster states become unstable and can therefore not be directly observed in experiments.

In Section~\ref{sec:Model}, we introduce the model of delay-coupled oscillators in both a general and experiment-specific framework, where the electrochemical oscillators are modelled by Stuart-Landau oscillators. We also investigate the phase space and parameter symmetries of the model. The results of a systematic analysis of the transitions between cluster states in the relaxation regime are presented in Section~\ref{sec:Results}. Finally, we discuss the results, how they relate to dynamics in the smooth regime, and further open questions in Section~\ref{sec:Discussion}.

\section{\label{sec:Model}Model}
In this section, we will present the model behind the numerical analysis. First, we will review the general framework of coupled Stuart-Landau oscillators in the context of cluster synchronization. In Section~\ref{subsec:coupling} we generalize the coupling scheme in the model for cases of non-sinusoidal coupling. Then, in Section~\ref{subsec:symmetries}, we describe clustering solutions from the model and discuss their symmetries.

\subsection{\label{subsec:SL}General framework}
Let us consider $N$ supercritical Hopf normal forms, also known as \textit{Stuart-Landau oscillators}, that are unidirectionally coupled in a  ring configuration with transmission delays:
\begin{align}
  \label{eq:SL}
  \dot{z}_j(t) =& \left[\lambda + i \omega - (1 + i \gamma) \left|z_j(t)^2\right| \right] z_j(t) \\ \nonumber
                &+ K\,z_{(j+1)\mod{N}}(t-\tau),
\end{align}
where $z_j(t)=r_j(t) e^{i\varphi_j(t)}\in\mathbb{C}$ denotes the $j$th complex dynamical variable, $j=1,\dots,N$. $\lambda$, $\omega\neq 0$, and $\gamma$ are real-valued constants, where $\gamma$ couples the frequency to the  oscillation amplitude. The parameters $K\in\mathbb{R}$ and $\tau$ denote the coupling strength and delay, respectively. In this study, we will focus on the formation of synchronized behavior as the delay is varied. In polar coordinates, that is, using radius $r_j$ and phase $\varphi_j$ variables, Eq.~\eqref{eq:SL} can be rewritten as follows:
\begin{subequations}
\label{eq:SL_polar}
\begin{align}
  \dot{r}_j(t) =& \left[\lambda - r_j(t)^2\right]r_j(t) \\\nonumber
                &+ K r_{j+1}(t-\tau) \cos\left[\varphi_{j+1}(t-\tau)-\varphi_j(t)\right],\\
  \dot{\varphi}_j(t) =& \omega - \gamma r_j(t)^2 \\\nonumber
                &+ K \frac{r_{j+1}(t-\tau)}{r_j(t)} \sin\left[\varphi_{j+1}(t-\tau)-\varphi_j(t)\right],
\end{align}
\end{subequations}
where all indexes have to be taken modulo $N$ throughout this paper.

\subsection{\label{subsec:coupling}Experiment-driven coupling}

As shown in Ref.~[\onlinecite{BLA13}], this simple model reproduces the correct clustering patterns and transitions between different patterns observed in the experiments for the smooth regime very well. However, in contrast to the smooth regime, the relaxation regime is not near the Hopf bifurcation that initiates oscillations. The coupling between oscillators can, therefore, not be assumed to be sinusoidal and we require a more general representation of the coupling. We rewrite the model using general interaction functions $H_i$, $i\in\{r, \varphi\}$:
\begin{subequations}
\label{eq:SL_polar_general}
\begin{align}
  \dot{r}_j(t) =& \left[\lambda - r_j(t)^2\right]r_j(t) \\\nonumber
                &+ K r_{j+1}(t-\tau)\, H_r\left[\varphi_{j+1}(t-\tau)-\varphi_j(t)\right],\\
  \dot{\varphi}_j(t) =& \omega - \gamma r_j(t)^2 \\\nonumber
                &+ K \frac{r_{j+1}(t-\tau)}{r_j(t)}\, H_\varphi\left[\varphi_{j+1}(t-\tau)-\varphi_j(t)\right].
\end{align}
\end{subequations}
In this paper, we use the radial and angular interaction functions as determined in Ref.~[\onlinecite{BLA13}] (cf. Section~V therein). In short, the experimentally obtained functions $H_r$ and $H_\varphi$ can each be approximated by a fifth-order Fourier series, which yields the following system of equations:
\begin{widetext}
\begin{subequations}
\label{eq:SL_polar_Fourier}
\begin{align}
  \dot{r}_j(t) =& \left[\lambda - r_j(t)^2\right]r_j(t) 
        + K r_{j+1}(t-\tau)
  % \sum_{n=1}^{4}a_{jn} r_n(t-\tau)
        \left(\sum_{n=0}^{5}
        a_{n,r} \cos\left\{n\left[\varphi_{j+1}(t-\tau)-\varphi_j(t)\right]\right\}
        + b_{n,r} \sin\left\{n\left[\varphi_{j+1}(t-\tau)-\varphi_j(t)\right]\right\}
    \right),\\
  \dot{\varphi}_j(t) =& \omega_j - \gamma r_j(t)^2
    + K \frac{r_{j+1}(t-\tau)}{r_j(t)} 
        \left(\sum_{n=0}^{5}
        a_{n,\varphi} \cos\left\{n\left[\varphi_{j+1}(t-\tau)-\varphi_j(t)\right]\right\}
        + b_{n,\varphi} \sin\left\{n\left[\varphi_{j+1}(t-\tau)-\varphi_j(t)\right]\right\}
    \right),    
\end{align}
\end{subequations}
\end{widetext}
where all indices have to be taken modulo $N$. 
The Fourier coefficients $a_{n,r}$, $b_{n,r}$, $a_{n,\varphi}$, and $b_{n,\varphi}$ that best match experimentally observed interaction functions, can be found in Appendix~\ref{app:Fourier}. The values of the other parameters in the model are chosen to match the experiments, as discussed in Ref.~[\onlinecite{BLA13}].

\subsection{\label{subsec:symmetries}Clustering patterns and their symmetries}

We investigate the occurrence of cluster states that are characterized by a common collective frequency $\Omega_m$. The cluster states come in two flavors.
\textit{Primary states} have a collective amplitude and equal phase lags $\Delta\varphi_m = 2\pi m/N$ between neighboring oscillators, i.e., $r_j\equiv r_{0,m}$ and $\varphi_j(t)=\Omega_m t + j \Delta\varphi_m$. 
\textit{Secondary states} differ in amplitudes and phase lags, but still oscillate with a common frequency.

The integer $m=0,\dots,N-1$ labels the specific primary cluster state.
Following Ref.~[\onlinecite{golubitsky12}], these states may be characterized as: discrete rotating waves when $m$ is coprime to $N$; discrete standing waves when $m=0$; and a discrete alternating wave when $m=N/2$.
In our case of $N=4$, $m=0$ corresponds to the in-phase state (cf. Fig.~\ref{fig:exp_primary}(a)), $m=1$ to the splay state (cf. Fig.~\ref{fig:exp_primary}(b)), $m=2$ to the 2-cluster state (cf. Fig.~\ref{fig:exp_primary}(c)), and $m=3$ to the reverse splay state (cf. Fig.~\ref{fig:exp_primary}(d)). 

Let us briefly discuss the symmetries present in the dynamical system \eqref{eq:SL_polar_Fourier} for the case $N=4$.  
As we will see throughout our analysis, it is helpful to relate the different clustering patterns observed in the experiments to their symmetry groups. This will aid us in our goal of understanding how they are related to each other and how the system transitions from one state to another.

There are two phase space symmetries, namely the symmetry under a global shift of all phase variables by a common phase $\varphi_s$, as well as the symmetry under cyclic permutations of the coordinate indices. More precisely, if we are given a particular solution $(\mathbf{r}(t),\pmb{\varphi}(t))^T$, with $\mathbf{r}(t)=(r_1(t),\ldots,r_4(t))^T$ and $\pmb{\varphi}(t)=(\varphi_1(t),\ldots,\varphi_4(t))^T$, then the functions defined via
\begin{equation}
\label{eq:Z4symm}
\begin{pmatrix}
\hat{\mathbf{r}}_{(\varphi_s,j)}(t) \\
\hat{\pmb{\varphi}}_{(\varphi_s,j)}(t)
\end{pmatrix} = 
\begin{pmatrix}
G^j\mathbf{r}(t) \\
G^j{\pmb{\varphi}}(t) + \varphi_s \mathbf{b}
\end{pmatrix} 
\end{equation}
are also solutions of \eqref{eq:SL_polar_Fourier} for all $\varphi_s \in [0,2\pi)$ and $j\in \{0,\ldots,3\}$. Here $G$ is the unidirectional coupling matrix defined in Eq.~\eqref{eq:Gmat} which cyclically permutes coordinate indices $j\to j-1$ and $\mathbf{b} = (1,1,1,1)^T$.  A pair of the form $(\phi_s,j)$ therefore performs a symmetry operation on the set of all solutions of \eqref{eq:SL_polar_Fourier}. The set of all symmetry operations $(\varphi_s,j)$ forms a group with neutral element $(0,0)$ and composition rule 
\begin{equation}
(\varphi_{s,1},j_1) \cdot (\varphi_{s,2},j_2) 
= (\varphi_{s,1}+\varphi_{s,2}\; \mathrm{mod}\; 2\pi,j_1 +j_2\; \mathrm{mod}\; 4).
\end{equation}
Let us denote this group by $H$. We note that $H$ is isomorphic to the direct sum $U(1) \oplus \mathbb{Z}_4$, where $U(1)=\{e^{i\phi}: \phi \in [0,2\pi) \}$ is the unitary group in one dimension and $\mathbb{Z}_4 = \mathbb{Z}/4\mathbb{Z} =  \{0,1,2,3\}$ is the cyclic group with four elements.  In mathematical terms we can also express the above observations by saying that the dynamical system \eqref{eq:SL_polar_Fourier} is $H$-equivariant \cite{balanov2006applied}. 

Now that we have identified the phase space symmetry group of our system, we can use it to classify the different solutions of the system based on their symmetry. We are particularly interested in cluster states, where all oscillators share a common frequency. For example, the splay state (i.e. the primary state with $m=1$ introduced in Sect.~\ref{subsec:SL}) is found to be invariant under the action of the element $\left(-\pi/2, 1\right)\in H$ since $\hat{\varphi}_{j}=\varphi_{j+1}+\pi/2=\varphi_{j}$. More generally, a primary state of index $m$ is invariant under a subgroup 
\begin{equation}
H_{4,m} = \{(- mj \pi/2, j): j=0,\ldots,3\} \subset H 
\end{equation}
Note that the order of all $H_{4,m}$ is four and they are isomorphic to $\mathbb{Z}_4$.  The secondary states are less symmetric and the three options for invariant subgroups are:
\begin{equation}
H_{2,0} = \{ (0,0), (0,2) \}; H_{2,1} = \{ (0,0), (\pi,2) \};
H_{1,0} = \{ (0,0) \}.
\end{equation}
$H_{2,0}$ and $H_{2,1}$ are isomorphic to $\mathbb{Z}_2$ and $H_{1,0}$ is the trivial group with only the neutral element. 
The groups of the clustering patterns observed in the experiments are summarized in Table~\ref{tab:Groups}. Note that the pattern corresponding to $H_{2,1}$, which was not observed in the experiments, will be discussed below in Section~\ref{subsec:open2cl}.

\begin{table}[htb!]
    \centering
    \caption{Symmetry groups of cluster states.}
    \label{tab:Groups}
    \begin{tabular}{lcc}
    \hline
    \hline
    Pattern & Group & Representation\\
    \hline
        In-phase  & $\mathbb{Z}_4$ & $H_{4,0}$ \\
        Splay & $\mathbb{Z}_4$ & $H_{4,1}$ \\
        2-cluster  & $\mathbb{Z}_4$ & $H_{4,2}$ \\
        Reverse splay  & $\mathbb{Z}_4$ & $H_{4,3}$ \\
        Compressed 2-cluster  & $\mathbb{Z}_2$ & $H_{2,0}$ \\
        Compressed splay & $\mathbb{Z}_1$ & $H_{1,0}$ \\
        Compressed reserve splay & $\mathbb{Z}_1$ & $H_{1,0}$ \\
        Open 2-cluster & $\mathbb{Z}_1$ & $H_{1,0}$ \\
    \hline
    \end{tabular}
\end{table}

In addition to the phase space symmetry, the dynamical system \eqref{eq:SL_polar_Fourier} also possesses an important parameter symmetry for cluster states.  For a particular $\tau$ let us assume that we are given a solution of \eqref{eq:SL_polar_Fourier} of the form
\begin{equation}
    \begin{pmatrix}
        \mathbf{r}(t) \\
        \pmb{\varphi}(t)
    \end{pmatrix} = 
    \begin{pmatrix}
        \mathbf{r}_0 \\
        \pmb{\varphi}_0+ \Omega t \mathbf{b}
    \end{pmatrix}, 
\end{equation}
with $\mathbf{r}_0, \pmb{\varphi}_0 \in \mathbb{R}^4$ and $\Omega \in \mathbb{R}$. With $\mathbf{c}= (0,\pi/2, \pi, 3\pi /2)^T$ it then follows that 
\begin{equation}
    \begin{pmatrix}
        \hat{\mathbf{r}}(t) \\
        \hat{\pmb{\varphi}}(t)
    \end{pmatrix} = 
    \begin{pmatrix}
        \mathbf{r}_0 \\
        \pmb{\varphi}_0+ \mathbf{c} + \Omega t \mathbf{b}
    \end{pmatrix} 
\end{equation} 
is a solution of \eqref{eq:SL_polar_Fourier} for $\hat{\tau}= \tau + \pi / (2 \Omega)$. This can be seen by noting that 
\begin{align}
&\hat{\varphi}_{j+1}(t-\hat{\tau})- \hat{\varphi}_{j}(t)\\ 
&= 
\varphi_{0,j+1} + (j+1)\frac{\pi}{2} - (\varphi_{0,j} + j\frac{\pi}{2}) - \Omega \left(\tau + \frac{\pi}{2 \Omega}\right) \\
&= \varphi_{j+1}(t-\tau)- \varphi_{j}(t).
\end{align} 

It is important to note that because of the presence of the vector $\mathbf{c}$ in the transformation, the symmetry group of the cluster state changes. For example, an in-phase state at delay $\tau$ which is invariant under $H_{4,0}$ will become a splay state at $\hat{\tau}$ with symmetry $H_{4,1}$.  More generally, $H_{4,m}$ becomes $H_{4,m+1}$ and the groups $H_{2,0}$ and $H_{2,1}$ get exchanged.

Repeated application of the above transformation yields a parameter symmetry group which is isomorphic to $\mathbb{Z}$. Explicitly, for every $j\in \mathbb{Z}$ the transformation
\begin{equation}
\label{eq:parasym}
    \begin{pmatrix}
        \hat{\mathbf{r}}_0 \\
        \hat{\pmb{\varphi}}_0 \\
        \hat{\tau} \\
        \hat{\Omega}
    \end{pmatrix} = 
    \begin{pmatrix}
        \mathbf{r}_0 \\
        \pmb{\varphi}_0+ j \mathbf{c} \\
        \tau + j \frac{\pi}{2\Omega} \\
        \Omega
    \end{pmatrix} 
\end{equation}
describes a parameter transformation between cluster states in the $(\tau,\Omega)$ plane.  As our original system is only defined for positive $\tau$, we will restrict ourselves to that case.  We also note that stability is not necessarily conserved under this parameter symmetry. This is in contrast to phase space symmetry, where stability information \emph{is} preserved.

\section{\label{sec:Results}Analysis}

In the following, we explain the nature of the experimentally observed, secondary cluster states by exploring the transitions between different primary cluster states. We demonstrate these transitions by means of bifurcation diagrams for the model given by Eqs.~\eqref{eq:SL_polar_Fourier} with $N=4$. We will demonstrate that already this small number gives rise to rich dynamical scenarios. 

The results in this paper are obtained numerically with the continuation software \textsc{DDE-Biftool} \cite{engelborghs00,SIE14a}. Throughout this section, we fix the system parameters as $\lambda=2.89$, $\omega=2.43$ and $K=0.189$, as in the previous work on the same experimental set-up\cite{BLA13}.
In order to facilitate a thorough analysis, we reduce the model by replacing the phase variables with phase-difference variables (see Appendix~(\ref{app:reduction}) for details). In particular, phase differences with respect to oscillator 4 are considered, so that oscillator 4 becomes the reduced system's rotating frame of reference.
This means that periodic solutions are represented by relative equilibria; for example, a 2-cluster periodic solution is represented by a relative equilibrium, where all phase differences are $\pi$. Therefore, in the bifurcation diagrams below, bifurcations of period solutions are represented by their counterpart bifurcations of relative equilibria, e.g., torus bifurcations are represented by Hopf bifurcations.

\begin{figure}[t]
\includegraphics[width=\linewidth]{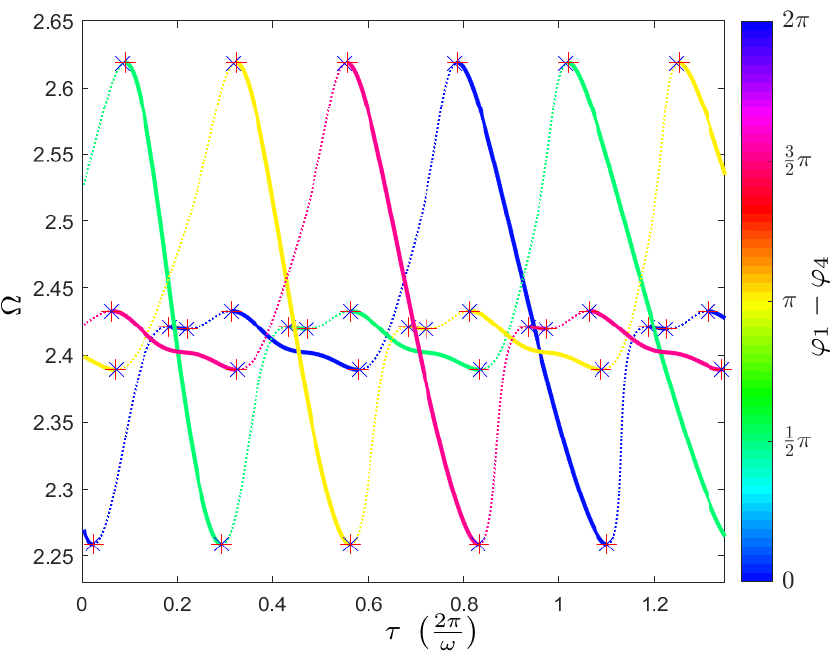}
\caption{Bifurcation diagram of primary states in the $(\tau,\Omega)$-plane. Solid/dotted curves represent stable/unstable solutions. Curve color indicates the phase difference between oscillators 1 and 4. Red plus signs/blue crosses denote Hopf/pitch-fork bifurcations. Parameters are $\lambda=2.89$, $\omega=2.43$ and $K=0.189$.}
\label{fig:bif_primary}
\end{figure}

Figure~\ref{fig:bif_primary} shows the bifurcation diagram of primary states in the $(\tau,\Omega)$-plane, that is, we vary the coupling delay and depict the collective frequency of the synchronized cluster state. Red plus signs and blue crosses denote Hopf and pitch-fork bifurcations, respectively, and stability/instability of the states are indicated by solid/dotted line styles. 
For better visualization, we use a color code to highlight the phase difference between oscillators 1 and 4, that is, $\varphi_1-\varphi_4 \in [0,2\pi]$.
Therefore, the blue, green, yellow, and red curves correspond to in-phase, splay, 2-cluster, and reverse splay states, respectively.  The different states are connected via the parameter symmetry \eqref{eq:parasym}, and this explains the striking similarity between the various branches. 

One can see how the different cluster states change in their collective frequency and how the stable branches that are experimentally accessible (see Ref.~[\onlinecite{BLA13}] for full branches) are connected by unstable branches. Stable branches connect the maximum and minimum of the collective frequency curves and also appear at intermediate frequencies around the intrinsic frequency $\omega$. The transitions between stability and instability occur via pairs of Hopf and pitchfork bifurcations in very close proximity to each other. We will use Fig.~\ref{fig:bif_primary} as an overall reference plot and explore how secondary states emerge in that bifurcation diagram. For this purpose, we will display the curves of the respective secondary states on top of the relevant primary state and adjust the range of the $\tau$-axis accordingly.

As a brief guide, we consider the secondary states: compressed 2-cluster state (Sec.~\ref{subsec:comp2cl}), compressed reverse splay state (Sec.~\ref{subsec:comprsp}), compressed splay state (Sec.~\ref{subsec:compsp}), and open 2-cluster state (Sec.~\ref{subsec:open2cl}). Aiming for an intuitive understanding, we also provide examples of solutions along the branches in the phase plane, which depict steps along the transition between different primary states.

\subsection{\label{subsec:comp2cl}Compressed 2-cluster state}

\begin{figure}[t]
\includegraphics[width=\linewidth]{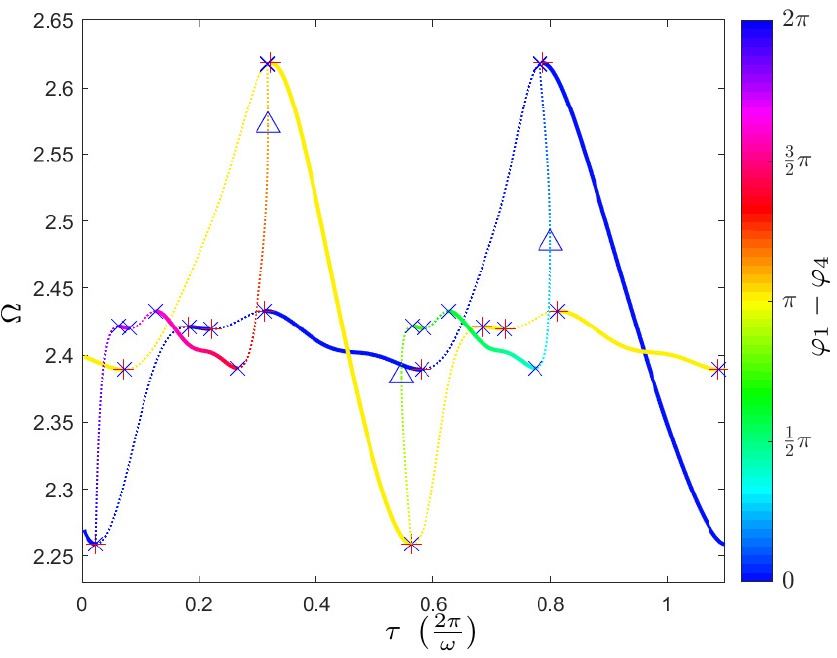}
\caption{Bifurcation diagram showing transitions between the 2-cluster and in-phase states in the $(\tau,\Omega)$-plane. Solid/dotted curves represent stable/unstable solutions. Curve color indicates the phase difference between oscillators 1 and 4. Plus signs, crosses, and triangles denote Hopf, pitch-fork, and fold bifurcations, respectively. Other parameters as in Fig.~\ref{fig:bif_primary}.
}
\label{fig:bif_inp_2cl}
\end{figure}

Let us start by considering the compressed 2-cluster state, shown schematically in Fig.~\ref{fig:exp_secondary}(b). This is a configuration, where oscillators group in alternating pairs \{1,3\} and \{2,4\}. The compression arises from a phase lag between the two pairs that is different from $\pi$.

\begin{figure}[t]
\includegraphics[width=\linewidth]{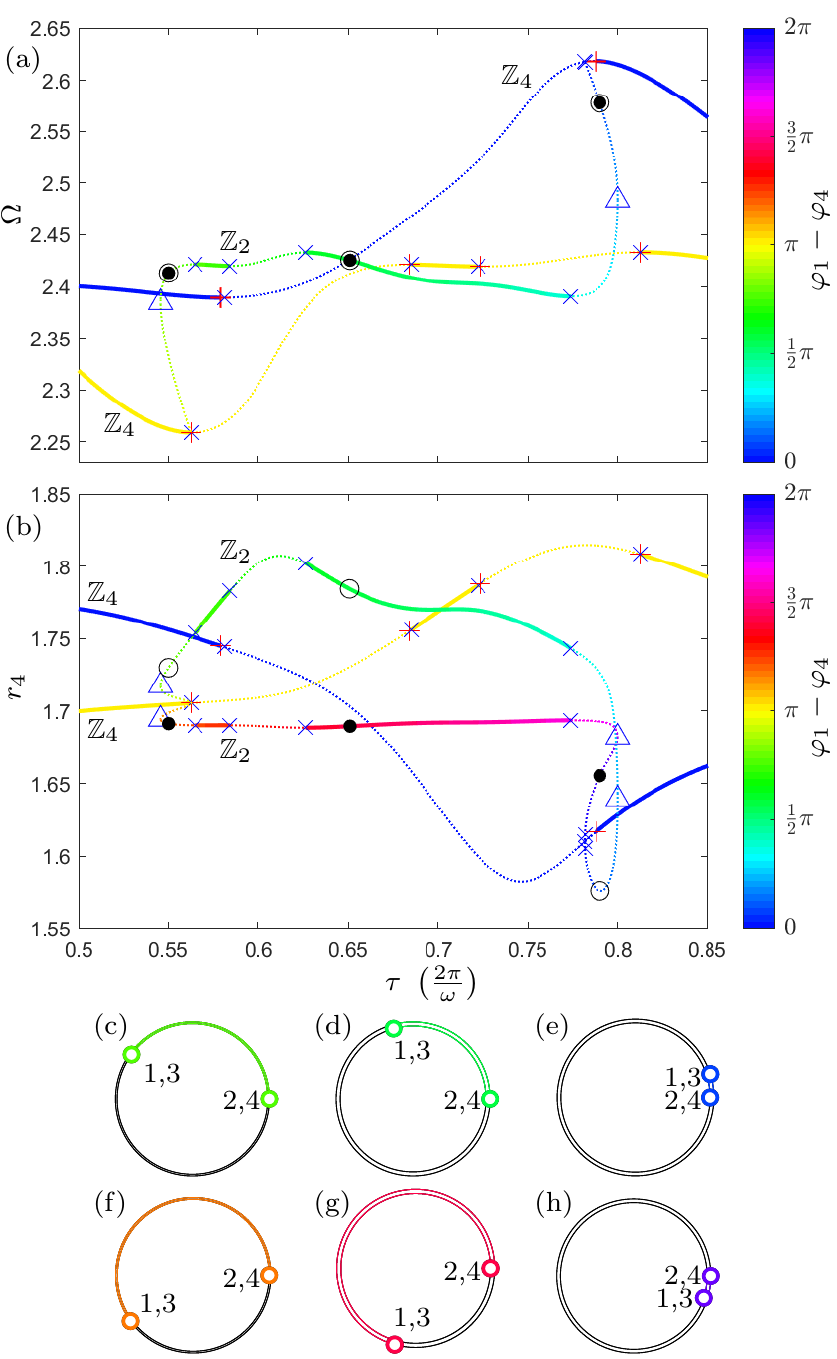}
\caption{Bifurcation diagram showing a transition between the 2-cluster and in-phase state in (a) the 
 $(\tau,\Omega)$-plane and (b) the $(\tau,r_4)$-plane. Solid/dotted curves represent stable/unstable solutions. Curve color indicates the phase difference between oscillators 1 and 4. Plus signs, crosses, and triangles denote Hopf, pitch-fork, and fold bifurcations, respectively. Open and closed circles denote the solutions shown in the phase plane in panels (c)--(e) and (f)--(h), respectively. The colored circular arcs highlight the phase differences between oscillators 1 and 4. Other parameters as in Fig.~\ref{fig:bif_primary}.
 }
\label{fig:bif_comp2cl}
\end{figure}

Figure~\ref{fig:bif_inp_2cl} illustrates how compressed 2-cluster states emerge as connecting branches between the in-phase (blue) and 2-cluster solutions (yellow). The bifurcation diagram highlights the recursive nature of these connecting branches: every time the in-phase solutions encounter a pitch-fork bifurcation, a branch of compressed 2-cluster solutions is born, which terminates at a pitch-fork bifurcation of 2-cluster solutions. This diagram also highlights how the phase differences along the primary branches are fixed. For example, the phase difference $\varphi_1 - \varphi_4$ is either fixed at $0$ or $\pi$. On the other hand, along the secondary branches the phase differences vary; providing the first clue to their role as transitional states. Here, the phase difference $\varphi_1 - \varphi_4$ varies continuously between $0$ and $\pi$ or $\pi$ and $2\pi$.

Figure~\ref{fig:bif_comp2cl} provides a closer look at the transition between in-phase and 2-cluster states. Panel~(a) is a zoomed-in version of the bifurcation diagram in the $(\tau,\Omega)$-plane shown in Fig.~\ref{fig:bif_inp_2cl}. Here, we see clearly that the branch of compressed 2-cluster states undergoes pitchfork bifurcations, creating intervals of stable solutions, which can be observed experimentally. Indeed, the stable interval of compressed 2-cluster states with approximately $\tau\in[0.63,0.76]$ was shown experimentally in Ref.~[\onlinecite{BLA13}] (see Fig.~8a therein). However, it is not immediately clear from panel~(a) that this branch emerges from pitchfork bifurcations, and that there are in fact two branches of compressed 2-cluster states that overlap each other in the $(\tau,\Omega)$-plane. Therefore, in panel~(b), we show the same diagram in the $(\tau,r_4)$-plane. This allows us to distinguish the different states via their respective radius. Now, we see the two branches of compressed 2-cluster states (with $\mathbb{Z}_2$-symmetry) that connect the in-phase and 2-cluster branches (with $\mathbb{Z}_4$-symmetry). Considering the color scheme, we see that one transition involves the phase difference $\varphi_1-\varphi_4$ increasing from $\pi$ to $2\pi$, while the other involves the same phase difference decreasing from $\pi$ to $0$. 

Panels~(c)--(e) of Fig.~\ref{fig:bif_comp2cl} are example solutions, corresponding to the open circles in the bifurcation diagrams, and panels~(f)--(h) are the symmetrically-related solutions, corresponding to the closed circles. Panels~(c) and~(f) both appear similar to a 2-cluster state, except that now the $\varphi_{j+1}-\varphi_j=\varphi_j-\varphi_{j-1}$ symmetry is broken as a result of the pitchfork bifurcation. Panels~(d) and~(g) are examples where the two clusters of oscillator pairs have moved closer to each other and are now stable. This example is typical of solutions observed experimentally, and corresponds to the schematic diagram shown in Fig.~\ref{fig:exp_secondary}(b); note the quantitative agreement of the delay times with the experimental observation. Finally, in panels~(e) and~(h), the two clusters are so compressed that they are almost in-phase.

In the context of symmetry, the compressed 2-cluster states are invariant under the group $H_{2,0}$. This is a subgroup of both $H_{4,0}$ and $H_{4,2}$, which are the symmetry groups of the in-phase state, and the 2-cluster state, respectively. It therefore follows that the compressed 2-cluster state can connect to those two states through pitchfork bifurcations.

\subsection{\label{subsec:comprsp}Compressed reverse splay state}

\begin{figure}[t]
\includegraphics[width=\linewidth]{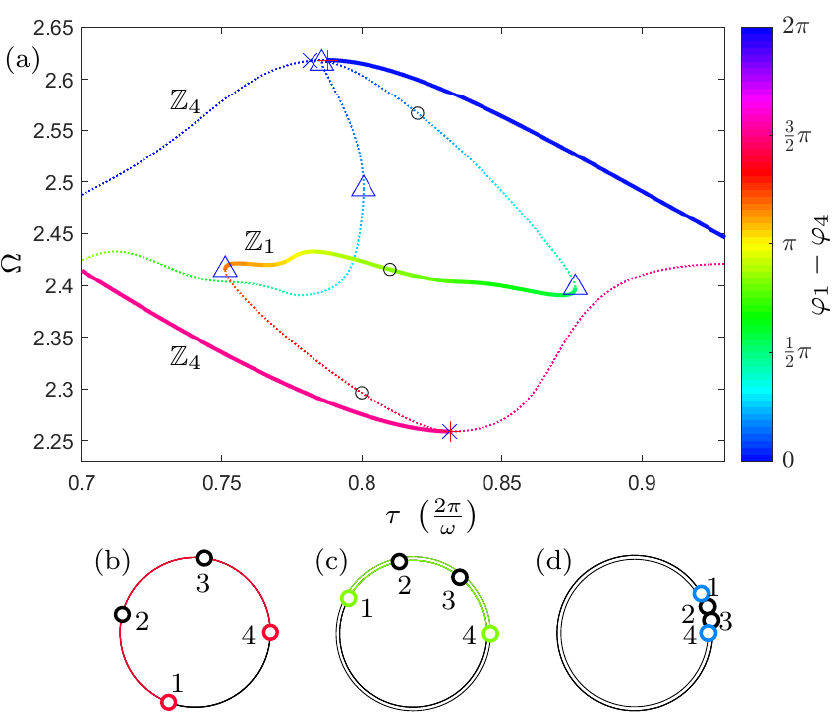}
\caption{Bifurcation diagram showing a transition between the reverse splay and in-phase state in the 
 (a) $(\tau,\Omega)$-plane. Solid/dotted curves represent stable/unstable solutions. Curve color indicates the phase difference between oscillators 1 and 4. Plus signs, crosses, and triangles denote Hopf, pitch-fork, and fold bifurcations, respectively. Circles denote the solutions shown in the phase plane in panels (b)--(d). The colored circular arcs highlight the phase differences between oscillators 1 and 4. Other parameters as in Fig.~\ref{fig:bif_primary}.
}
\label{fig:bif_comprsp}
\end{figure}

Figure~\ref{fig:bif_comprsp}(a) shows a branch of compressed reverse splay states, together with branches of reverse splay and in-phase states. The example solutions in panels~(b)--(d) depict a gradual transition from a reverse splay to an in-phase configuration. The stable solution shown in panel~(c) is typical of experimental observations and corresponds to the schematic diagram shown in Fig.~\ref{fig:exp_secondary}(d); again, note the quantitative agreement of $\tau$ with the experiment.

\begin{figure}[t]
\includegraphics[width=\linewidth]{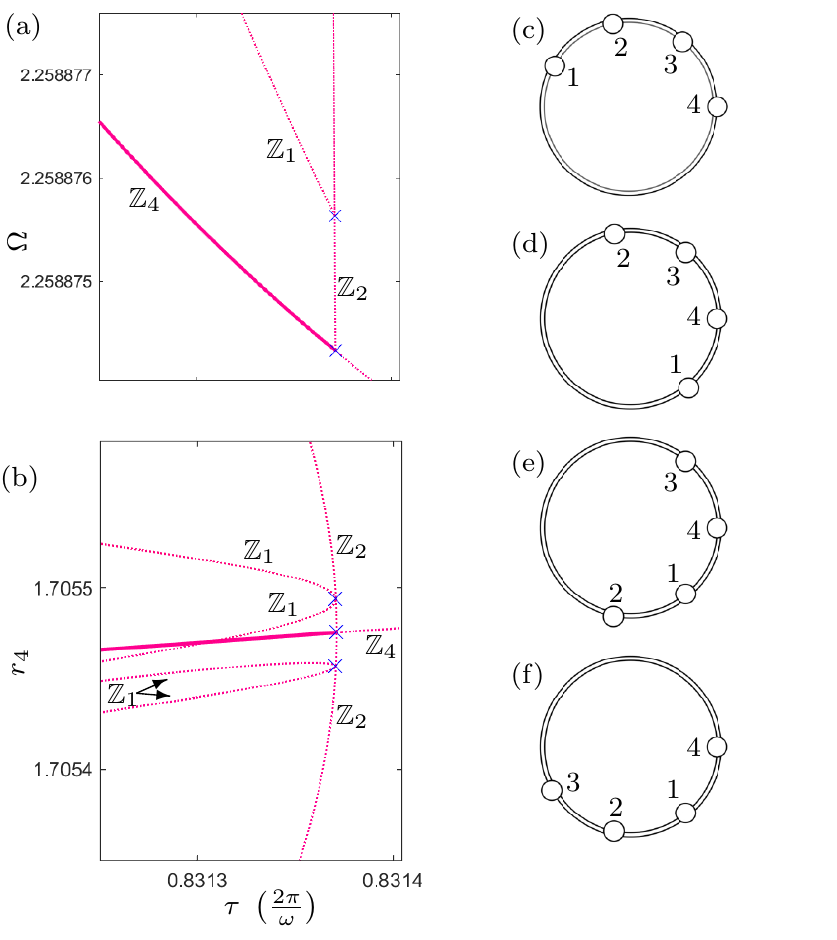}
\caption{Pitchfork bifurcations, denoted by blue crosses, of the reverse splay state in (a) the 
 $(\tau,\Omega)$-plane and (b) the $(\tau,r_4)$-plane. Solid/dotted curves represent stable/unstable solutions. Panels (c)--(f) are solutions resulting from the pitchfork bifurcations for the same parameters and collective frequency $\Omega$ as Fig.~\ref{fig:bif_comprsp}(c). Other parameters as in Fig.~\ref{fig:bif_primary}.
}
\label{fig:bif_pitchforks}
\end{figure}

In this case, however, the branch of secondary cluster states is not connected to the primary states branches by a simple pair of pitchfork bifurcations. This is clear from the fact that the primary branches are $\mathbb{Z}_4$, while this branch of secondary cluster states possesses only $\mathbb{Z}_1$ symmetry. Figure~\ref{fig:bif_pitchforks} provides further details on how the branch of compressed reverse splay states is related to the branch of reverse splay states. Panels~(a) and (b) show the reverse splay branch losing stability at a pitchfork bifurcation. At this initial pitchfork bifurcation, the $\varphi_{j+1}-\varphi_j=\varphi_j-\varphi_{j-1}$ symmetry is lost, and two branches of solutions with $\mathbb{Z}_2$ symmetry are born. Then, very shortly afterwards, further pitchfork bifurcations take place, at which the $\varphi_{j+2}-\varphi_j=\varphi_j-\varphi_{j-2}$ symmetry is lost, resulting in four branches with $\mathbb{Z}_1$ symmetry. This is in agreement with the \emph{equivariant branching lemma} \cite{golubitsky12}, which characterizes the isotropy subgroups of bifurcating solution branches. Only with both sets of symmetries gone, can all oscillators approach each other to form the compressed reverse splay state. As a result of the sets of pitchfork bifurcations, the compressed reverse splay branch shown in Fig.~\ref{fig:bif_comprsp}(a) is actually four symmetrically-related curves overlapping in the $(\tau,\Omega)$-plane. Figures~\ref{fig:bif_pitchforks}(c)--(f) show an example of four symmetrically-related solutions in the phase plane, each with the same parameters and collective frequency $\Omega$ as Fig.~\ref{fig:bif_comprsp}(c).

\begin{figure}[t]
\includegraphics[width=\linewidth]{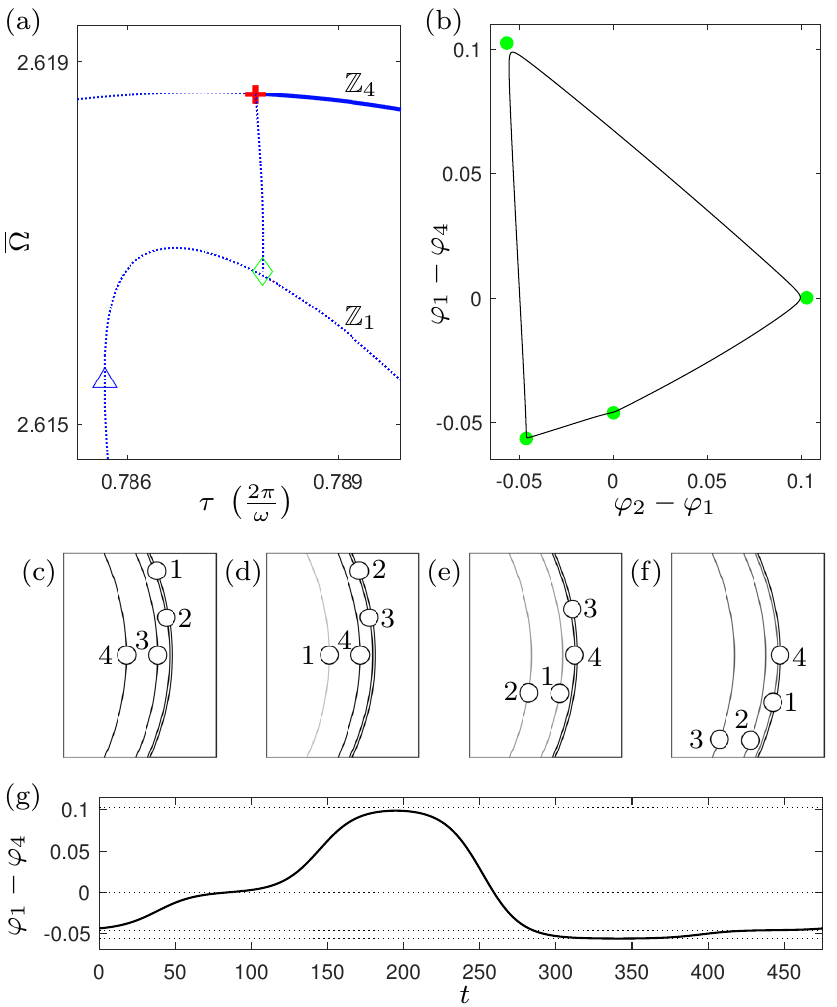}
\caption{(a) Bifurcation diagram showing the connection between the in-phase branch and secondary cluster states involving periodic solutions in the $(\tau,\overline\Omega)$-plane, where $\overline\Omega$ is mean collective frequency. 
Solid/dotted curves represent stable/unstable solutions. The cross, triangle, and diamond denotes a Hopf, fold, and heteroclinic bifurcation, respectively. (b) A periodic solution close to the heteroclinic bifurcation projected onto the $(\varphi_2-\varphi_4,\varphi_1-\varphi_4)$-plane, where green circles represent saddles. Panels~(c)--(f) show close-ups of the four saddles shown in panel~(b) in the $(x,y)$ phase plane. (g) Time series of $\varphi_1-\varphi_4$ corresponding to the solution shown in panel~(b). Dotted lines correspond to the four saddles. 
Other parameters as in Fig.~\ref{fig:bif_primary}.
}
\label{fig:bif_hopf_hetero}
\end{figure}

The connection between the compressed reverse splay branch and the in-phase branch in Fig.~\ref{fig:bif_comprsp}(a) also requires a more detailed inspection. Figure~\ref{fig:bif_hopf_hetero}(a) provides a zoomed-in view of the bifurcation diagram. Here, the vertical axis represents the mean collective frequency of the four oscillators, $\overline\Omega$, because the bifurcation diagram now involves periodic solutions. At the Hopf bifurcation, denoted by the red cross, the in-phase solutions lose stability and a branch of unstable periodic solutions emerge. To emphasize, these are periodic solutions of the phase \emph{differences} of the oscillators. In the original non-rotating frame of reference used in Eqs.~\eqref{eq:SL_polar_Fourier}, these solutions correspond to solutions on a torus. 
The periodic solutions, born at the Hopf bifurcation, fulfill the $\mathbb{Z}_4$ symmetry condition $\varphi_{j+1}(t)=\varphi_{j}(t+T/4)$ where $T$ is the period. This condition persists as they approach the four saddles on the four symmetrically-related $\mathbb{Z}_1$ branches of compressed reverse splay states, where a heteroclinic bifurcation takes place, denoted by the green diamond. A periodic solution approaching the heteroclinic bifurcation is shown in panel~(b) as a projection onto the $(\varphi_1-\varphi_4,\varphi_2-\varphi_1)$-plane, where green circles are the saddles. Panels~(c)--(f) show zoomed-in plots of the saddles in the phase plane, providing a clear illustration of the symmetry. Finally, we see in panel~(g) that the time series appears to plateau at the values of $\varphi_1-\varphi_4$ that correspond to the four saddles. This shows that, as expected with increasing period towards a heteroclinic bifurcation, the periodic solution begins to spend more time in the vicinity of the saddles.

\subsection{\label{subsec:compsp}Compressed splay state}

\begin{figure}[t]
\includegraphics[width=\linewidth]{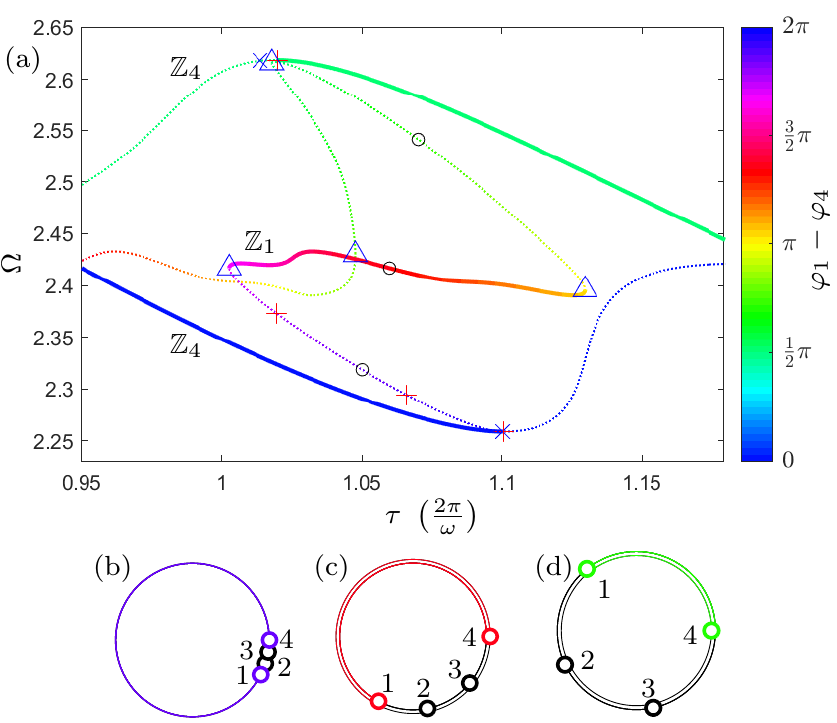}
\caption{Bifurcation diagram showing a transition between the in-phase and splay state in the 
 (a) $(\tau,\Omega)$-plane. Solid/dotted curves represent stable/unstable solutions. Curve color indicates the phase difference between oscillators 1 and 4. Plus signs and crosses, and triangles denote Hopf, pitch-fork, and fold bifurcations, respectively. Circles denote the solutions shown in the phase plane in panels (b)--(d). The colored circular arcs highlight the phase differences between oscillators 1 and 4. Other parameters as in Fig.~\ref{fig:bif_primary}.
}
\label{fig:bif_compsp}
\end{figure}

The compressed splay state arises in an analogous fashion to the compressed reverse splay state, except that here the transition is between the in-phase and splay primary branches. Figure~\ref{fig:bif_compsp}(a) shows the corresponding bifurcation diagram in the $(\tau,\Omega)$-plane, with a branch of compressed splay states connecting the two primary state branches.
In this case, the connection to the in-phase branch is via multiple pitchfork bifurcations, i.e., the mechanism demonstrated in Fig.~\ref{fig:bif_pitchforks}, while the connection to the splay branch is via a branch of periodic solutions, i.e., the mechanism demonstrated in Fig.~\ref{fig:bif_hopf_hetero}.

Panels~(b)--(d) show example solutions along the compressed splay branch, indicated by the circles in panel~(a). Again, the stable solution, shown in panel~(c), agrees with the schematic diagram of experimentally observed compressed splay states shown in Fig.\ref{fig:exp_secondary}(f) for the same value of $\tau$.  The striking similarity between Figs.~\ref{fig:bif_comprsp}(a) and~\ref{fig:bif_compsp}(a) is again explained through the connection of the branches via the parameter symmetry of \eqref{eq:parasym}. However, we stress that this connection only holds for cluster states and not limit cycles. It also does not preserve stability of states.  This explains that observed minor difference between Figs.~\ref{fig:bif_comprsp}(a) and~\ref{fig:bif_compsp}(a) which indicates that here the transition branch encounters a pair of Hopf bifurcations. 

\subsection{\label{subsec:open2cl}Open 2-cluster state}

\begin{figure}[t]
\includegraphics[width=\linewidth]{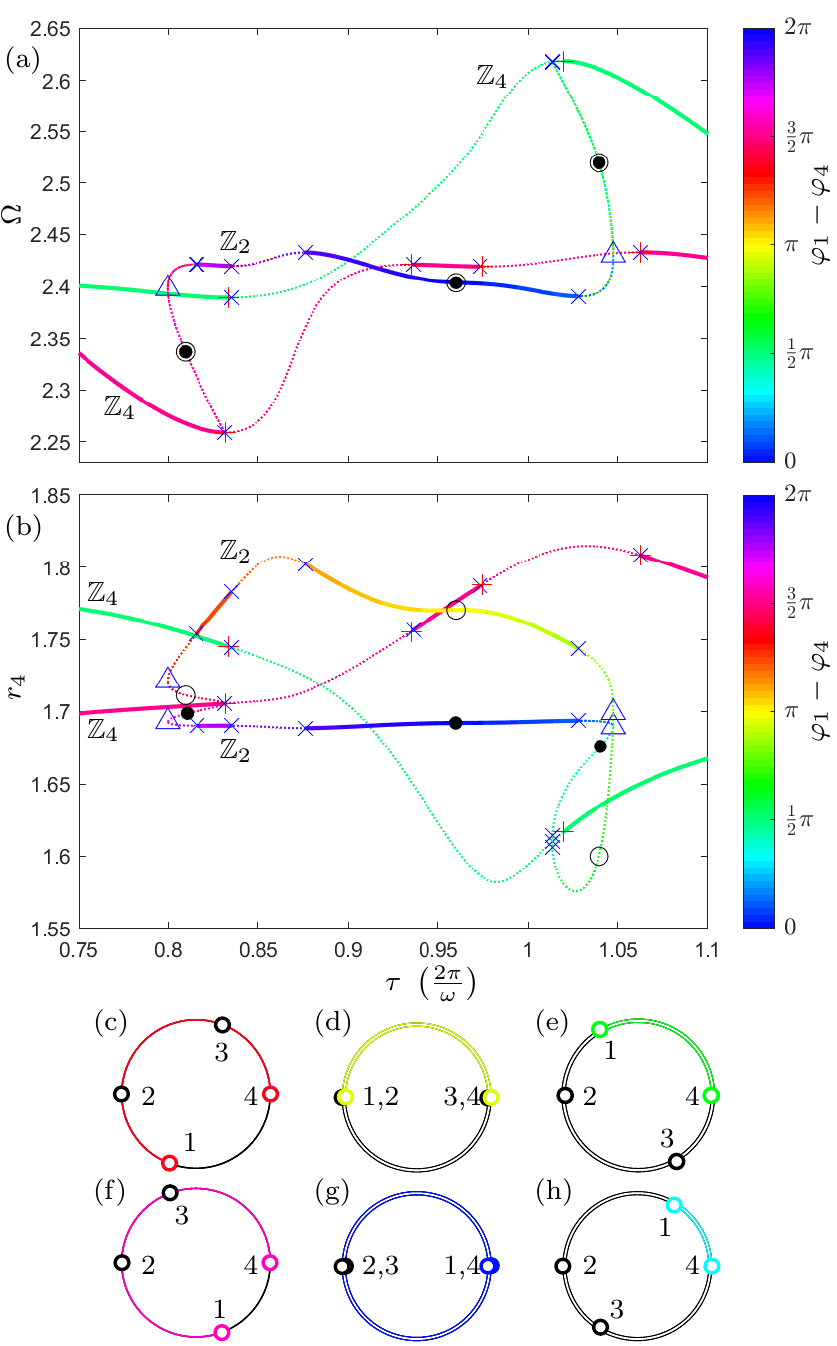}
\caption{Bifurcation diagram showing a transition between the reverse splay and splay state in (a) the 
 $(\tau,\Omega)$-plane and (b) the $(\tau,r_4)$-plane.  Solid/dotted curves represent stable/unstable solutions.} Curve color indicates the phase difference between oscillators 1 and 4. Plus signs, crosses, and triangles denote Hopf, pitch-fork, and fold bifurcations, respectively. Open and closed circles denote the solutions shown in the phase plane in panels (c)--(e) and (f)--(h), respectively. The colored circular arcs highlight the phase differences between oscillators 1 and 4. Other parameters as in Fig.~\ref{fig:bif_primary}.
\label{fig:bif_open2cl}
\end{figure}

In the model, the simplest way for open 2-cluster states to emerge is along branches connecting the splay and reverse splay branches. Figure~\ref{fig:bif_open2cl} shows two (symmetrically-related; $\mathbb{Z}_2$) open 2-cluster branches in (a) the $(\tau,\Omega)$-plane and (b) the $(\tau,r_4)$-plane that connects to the ($\mathbb{Z}_4$) primary branches at pitchfork bifurcations. Note that the curves of solutions take on the exact same forms as the curves in Fig.~\ref{fig:bif_comp2cl}, although different cluster types are involved. This is due to the parameter symmetry \eqref{eq:parasym}. While in Fig.~\ref{fig:bif_comp2cl} the connecting branch has $H_{2,0}$ symmetry, the branch in  Fig.~\ref{fig:bif_open2cl} has $H_{2,1}$ symmetry.  Therefore, throughout the transition from reverse splay to splay state in panels~(c)--(e) and~(f)--(h) of Fig.~\ref{fig:bif_open2cl} we see how oscillators 1 and 3 maintain a phase difference of $\pi$, as do oscillators 2 and 4. However, the phase differences between neighboring oscillators shift, allowing pairs of oscillators to approach each other. In panels~(d) and~(g) we see an example of pairs of neighboring oscillators with zero phase difference, resembling a case of a ``closed'' 2-cluster state. Note that this state is different to the primary 2-cluster state, where the clusters are given by pairs of oscillators $\{1,3\}$ and $\{2,4\}$ and all radii are equal.

\begin{figure}[t]
\includegraphics[width=\linewidth]{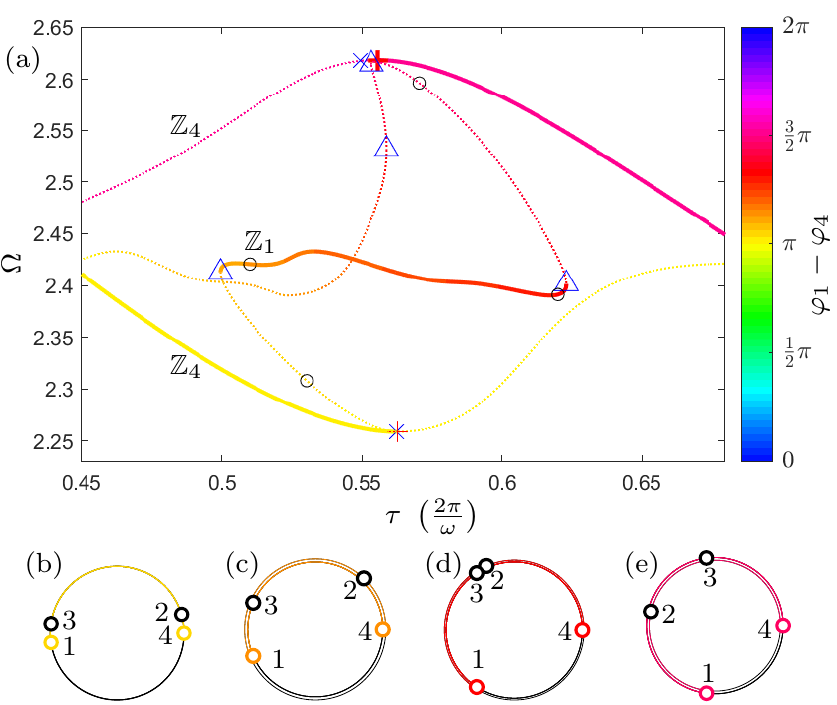}
\caption{Bifurcation diagram showing a transition between the 2-cluster and reverse splay state in the 
 $(\tau,\Omega)$-plane. Solid/dotted curves represent stable/unstable solutions. Curve color indicates the phase difference between oscillators 1 and 4. Plus signs, crosses, and triangles denote Hopf, pitch-fork, and fold bifurcations, respectively. Open circles denote the solutions shown in the phase plane in panels (b)--(e). The colored circular arcs highlight the phase differences between oscillators 1 and 4. Other parameters as in Fig.~\ref{fig:bif_primary}.
}
\label{fig:bif_open2cl_2}
\end{figure}

\begin{figure}[th]
\includegraphics[width=\linewidth]{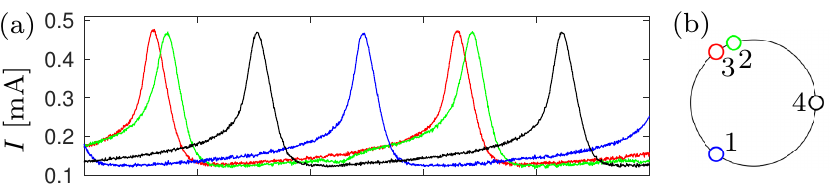}
\caption{Experimental time series (a) and schematic diagram (b) of a secondary cluster state in the relaxation regime: open 2-cluster state with $\tau= 0.62\times\left(\frac{2\pi}{\omega}\right)$. }
\label{fig:exp_secondary_open2cl}
\end{figure}

Although the secondary branch between the reverse splay and splay branches in Fig.~\ref{fig:bif_open2cl} provides a mechanism for solutions that could be described as open 2-cluster, they do not correspond to the open 2-cluster states observed in the electrochemical oscillator experiments. Comparison with the schematic diagram in Fig.~\ref{fig:exp_secondary}(h) reveals that the phase difference of $\pi$ between oscillators 1 and 3, as well as 2 and 4 is not present in the open 2-cluster states observed experimentally.

Figure~\ref{fig:bif_open2cl_2} shows the transition branch between the 2-cluster and reverse splay branches that provides the mechanism for the open 2-cluster states, as they appear in the experiments. The details of the connections between this open 2-cluster branch and the primary branches are analogous to Figs.~\ref{fig:bif_comprsp} and~\ref{fig:bif_compsp}. Generally, the phase differences between the paired oscillators in the 2-cluster state increase (see Fig.~\ref{fig:bif_open2cl_2}(b)) until they approach the reverse splay state (see Fig.~\ref{fig:bif_open2cl_2}(e)). Along the way, the clustering pattern appears as in the experiments (cf. panel~(c) with Fig.~\ref{fig:exp_secondary}(h)).

Comparing Fig.~\ref{fig:bif_open2cl_2}(c) and~(e) it becomes clear that along the branch oscillators 2 and 3 must pass each other. Panel~(d) shows an example of
a stable solution where oscillators 2 and 3 are very close together. Now, if we return to the experimental output and consider the open 2-cluster state just before it disappears as $\tau$ is being increased, we indeed find a cluster state (overlooked in Ref.~[\onlinecite{BLA13}], now shown in Fig.~\ref{fig:exp_secondary_open2cl}) that matches Fig.~\ref{fig:bif_open2cl_2}(d) and exists for the same value of $\tau$.

\section{\label{sec:Discussion}Discussion}

We have analyzed a mathematical model of delay-coupled oscillators to gain a deeper understanding of experimental results of a system of delayed-coupled electrochemical oscillators. Our analysis, conducted by numerical continuation, details the bifurcations involved in the appearance and disappearance of branches of various cluster types in experiments. Furthermore, we elucidate the role of the secondary cluster states as transitions between certain primary states. While primary states have symmetry groups  of order 4, secondary cluster states have symmetry groups of smaller order, and the connection between them is mediated through pitchfork bifurcations.

\begin{figure}[t]
\includegraphics[width=\linewidth]{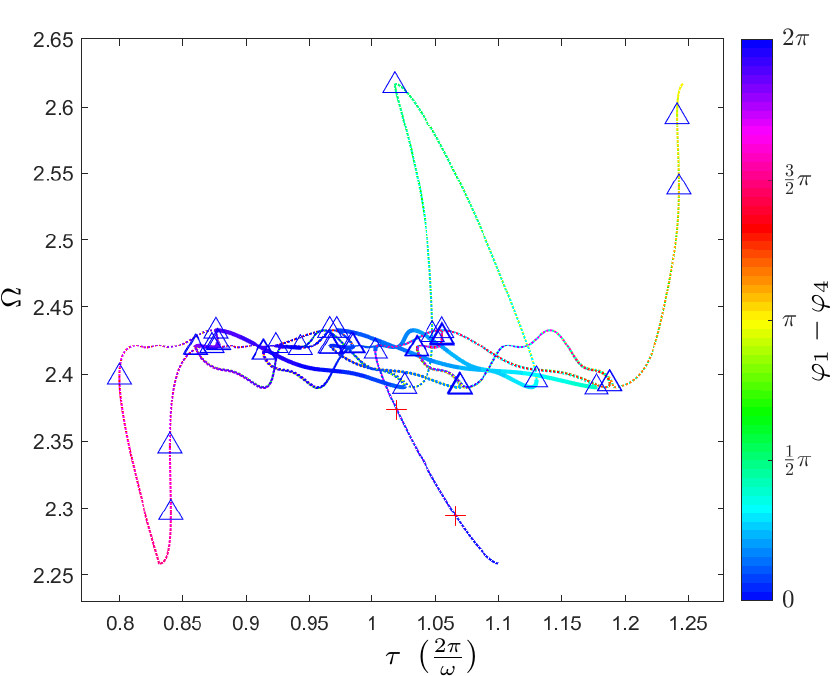}
\caption{The entire branch of secondary cluster states shown only partially in Fig.~\ref{fig:bif_compsp}. Solid/dotted curves represent stable/unstable solutions. Curve color indicates the phase difference between oscillators 1 and 4.} Plus signs and triangles denote Hopf and fold bifurcations, respectively. Other parameters as in Fig.~\ref{fig:bif_primary}.
\label{fig:bif_compsp_full}
\end{figure}

We explicitly demonstrate where each of the four types of experimentally-observed secondary clusters belong along the transitional routes between certain primary states.
As emphasized earlier, in each case of secondary state, what we show in the above figures are only the parts of the curves needed to explain the relevant transition. 
In fact, the curves are more complicated and intertwined. For example, Fig.~\ref{fig:bif_compsp_full} shows the complete branch with $\mathbb{Z}_1$ symmetry to which the compressed splay states in Fig.~\ref{fig:bif_compsp} belongs. Interestingly, it reveals that this single secondary branch actually provides continuous routes between all primary states via a plethora of fold bifurcations, since each extremum of the branch connects to one of the primary branches (not shown for clarity) by the same mechanisms demonstrated in Figs.~\ref{fig:bif_pitchforks} and~\ref{fig:bif_hopf_hetero}.

Our results have focused on the clustering patterns, as observed in the relaxation oscillation regime. But how do these results relate to the smooth oscillation regime where the oscillators first begin to oscillate upon increasing the potential $V_0$? Do the secondary states also exist in this regime? Why were the secondary states not observed there?

\begin{figure}[t]
\includegraphics[width=\linewidth]{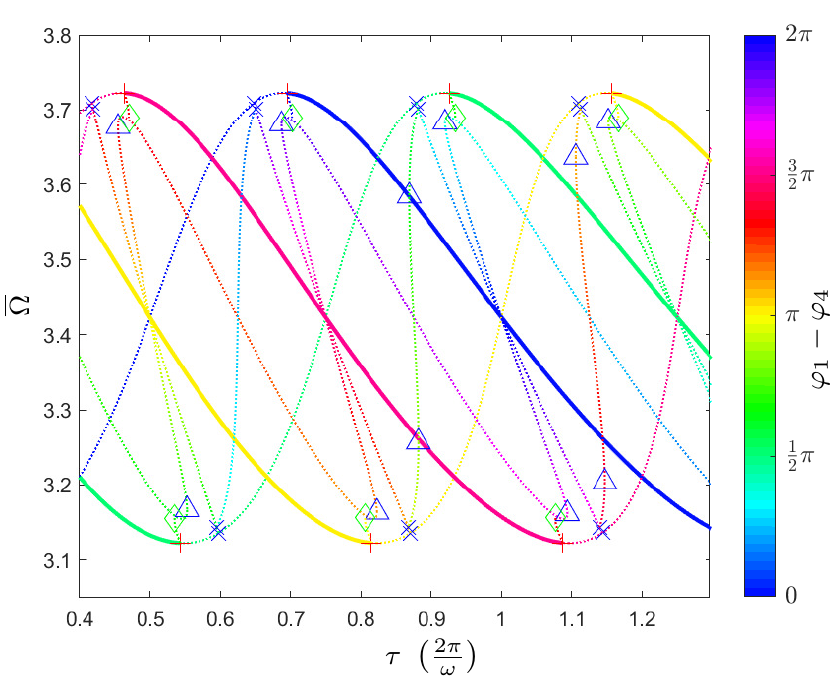}
\caption{Bifurcation diagram of smooth oscillators in the $(\overline\Omega,\tau)$-plane. Solid/dotted curves represent stable/unstable solutions. Curve color indicates the phase difference between oscillators 1 and 4. Plus signs, crosses, triangles, and diamonds denote Hopf, pitch-fork, fold, and heteroclinic bifurcations, respectively. Parameters: $\lambda=1.1025$, $\omega=3.4228$, and $K=0.3$.
}
\label{fig:bif_smooth}
\end{figure}

Figure~\ref{fig:bif_smooth} is the bifurcation diagram for the smooth oscillation regime in the $(\overline\Omega,\tau)$-plane. It reveals that, according to the model, the primary cluster states, observed in the experiments, do possess the same bifurcations as in the relaxation regime. Furthermore, similar to in our results above, we find branches of secondary clustering states that connect the different primary state branches via pitchfork bifurcations and small branches of periodic solutions. However, this interesting behavior of secondary states is always unstable. This provides an explanation for the lack of secondary states being observed in the experiments. In the smooth regime the secondary states are unstable and cannot be observed directly. It is only when the potential $V_0$ in the experimental setup is increased further, moving the system away from the smooth regime, that the bifurcation structure of the secondary state branches becomes more complex and creates pockets of stable solutions.

When we adapt the original model~(\ref{eq:SL_polar}) for general coupling forms in model~(\ref{eq:SL_polar_general}), we tailor the relevant terms to mimic the coupling observed in the experiments in the relaxation regime, as outlined in Section~\ref{subsec:coupling}. 
This turns out to be sufficient for achieving a very good agreement between the model and experiments. For example, in this paper we match examples of various cluster solutions from the model to experimental observations for the same value of $\tau$. Furthermore, the model reproduces the bifurcation structure inferred from the experiments with good quantitative agreement regarding oscillation frequency and $\tau$ (see Figs.~6 and~8 of Ref.~[\onlinecite{BLA13}]). 
Therefore, in contrast to the modelling of the coupling, the modelling of the local dynamics of the oscillators, which is not adapted in model~(\ref{eq:SL_polar_general}), appears to be less important for the study of the clustering behavior of the network. 

Further work could include bridging the gaps between the two snapshots we have of the smooth regime ($V_0=1.105$) and the relaxation regime ($V_0=1.2$). The key to this would be studying how the interaction functions evolve between the regimes. Another obvious avenue for further work would be to investigate the influence of the coupling strength $K$, which would undoubtedly be crucial for synchronization behavior.

The observation of the consistently close pairing of Hopf and pitch-fork bifurcations along the curves of primary states could warrant further study. Similar observations are made in other systems \cite{clerkin2014multistabilities,erzgraber2006compound}. How generic is this observation in relation to cluster synchronization? Is it even a dynamical requirement that a Hopf bifurcation occur in close proximity to the pitchfork bifurcations along the primary branches?

While we have found very good agreement between the secondary clustering patterns in the model and the experimental observations, there is one pattern in the model that was not observed experimentally; namely, the open 2-cluster state, as shown in Fig.~\ref{fig:bif_open2cl}(d) and (g). Of course, this could be due to the model not capturing some aspect of the network of coupled oscillators that exclude this particular state from observations. Alternatively, one possible explanation, which would require closer study, is that those solutions simply have relatively small basins of attraction compared to other co-existing solutions. For example, the range of $\tau$, for which the secondary branch in Fig.~\ref{fig:bif_open2cl} is stable, lies within the range of $\tau$ for which the in-phase branch is stable. Perhaps, when the oscillators transition between cluster states, some secondary branches are overlooked because of the size of their basin of attraction at the time of the transition.

Finally, we have remarked throughout the paper that the unstable solutions cannot be \emph{directly} observed in the experiments. One possible approach to indirectly observe these states experimentally would be to apply Pyragas control in order to stabilise the unstable solutions \cite{pyragas92}. This would require implementing a delayed self-feedback to each of the four oscillators and setting this delay time equal to the period of the unstable cluster state to be stabilized. Assuming the model to be accurate also for the unstable states, this period would be equal to $2\pi/\Omega$. One challenge might be finding an appropriate control strength to successfully stabilize the unstable solutions, nonetheless, it would be fasciniating to explore how well the model captures not only the stable, but also unstable clustering patterns of the electrochemical oscillations. 

\begin{acknowledgments}
The project was supported by University College Cork in the framework of the \textit{SEFS New Connections Grant Award} scheme.
PH acknowledges further support by \textit{Deutsche Forschungsgemeinschaft} (DFG) under project ID 434434223 - SFB 1461.
We are also grateful to two anonymous referees for valuable feedback, which helped us to improve the manuscript.
\end{acknowledgments}

\section*{Data Availability Statement}
The data that support the findings of this study are available from the corresponding author upon reasonable request.

% APPENDIX, IF NEEDED
%
\appendix

\section{\label{app:experiment}Experimental set-up}

The motivation for the theoretical study presented in this paper is to provide a thorough understanding of the transitions between different clustering patterns during experiments with electrochemical oscillators. The observed clustering patterns are introduced in Ref.~[\onlinecite{BLA13}], as well as numerous examples of transitions that occur as the control parameter, i.e. the coupling delay, is varied. Some additional (quantitatively different) runs are shown in this paper in Figs.~\ref{fig:exp_primary} and \ref{fig:exp_secondary}, and a (qualitatively different) run is shown in Fig.~\ref{fig:exp_secondary_open2cl}.
Therefore, for the convenience of the reader, we provide a brief summary of the experiment here. For further details, we refer the reader to Ref.~[\onlinecite{BLA13}].

We conduct experiments performed in an electrochemical cell consisting of four 1-mm-diameter Ni working electrodes (99.98\% pure), a Pt mesh counter electrode and Hg/Hg\textsubscript{2}SO\textsubscript{4}/K\textsubscript{2}SO\textsubscript{4} (sat) reference electrode with a 3M H\textsubscript{2}SO\textsubscript{4} electrolyte.  
We electrically couple the four electrodes in a unidirectional ring and measure the electrochemical dissolution currents  $I_j, j=1,2,3,4,$ via a zero resistance amperemeters (ZRAs). We select four oscillators with similar uncoupled frequencies from an array of 64 oscillators. We can vary the character of the oscillators with our choice of applied voltage. Lower voltages closer to the Hopf bifurcation are nearly harmonic; higher voltages produce higher harmonic oscillations that appear less sinusoidal. These two dynamical behaviors are the above-mentioned \textit{smooth} and \textit{relaxation} oscillations, respectively. In these experiments, we observe smooth and relaxation oscillators for an applied voltage of $V_0=1.105$~V and 1.2~V, respectively. 

Negligible intrinsic electrical interactions exist between uncoupled oscillators \cite{ZHA08}. The startup or shutdown of one oscillator does not alter the behavior of the others and oscillator dynamics are independent when uncoupled.

We introduce coupling of the form
\begin{equation}
    V_j(t)=V_0+\delta V_j(t),
\end{equation}
where $V_j,j=1,2,3,4$, denotes the voltage between working and reference electrodes, and $\delta V_j$ is the change in circuit potential of the \textit{j}th element due to feedback. The feedback voltages are
\begin{equation}
    \delta V_j(t)=K\displaystyle\sum_{n=1}^{4} g_{jn}[V_n(t-\tau)-R_P \hat{I}_n(t-\tau)],    
\end{equation}
where $R_P=650\,\Omega$ is a resistance, $K$ denotes the overall coupling gain, and $\tau$ is the coupling time delay;
we impose time delayed coupling via the real-time data acquisition system with a multichannel potentiostat. 
$\hat{I}_j$ is the normalized current measured by the ZRAs, such that
\begin{equation}
    \hat{I}_j(t)=\frac{\overline{A}}{A_j}(I_j(t)-\overline{I}_j),
\end{equation}
where $A_j$ and $\overline{I}_j$ are the amplitude and mean current of oscillator $j$ and $\overline{A}$ is the mean amplitude of the population, $\sum_{n=1}^{4}I_n^{\textrm{max}}/4$.

The coupling structure is determined by $g_{jn}$, which belongs to the adjacency matrix $\mathbf{G}$. We apply unidirectional coupling in a four-member ring with
\begin{equation}
\label{eq:Gmat}
\mathbf{G} =
\begin{pmatrix}
0 & 1 & 0 & 0\\
0 & 0 & 1 & 0\\
0 & 0 & 0 & 1\\
1 & 0 & 0 & 0
\end{pmatrix},
\end{equation}
which is implemented via the multichannel potentiostat.

\section{\label{app:Fourier} Fourier coefficients of the interaction functions}
The Fourier coefficients of the radial and angular interaction functions, which are used in Eqs.~\eqref{eq:SL_polar_Fourier}, are obtained in Ref.~[\onlinecite{BLA13}] and are summarized in Tab.~\ref{tab:Fourier} for the reader's convenience. They are scaled to normalize the angular interaction function $\max \left|H_\varphi\right| = 1$.

\begin{table}[htb!]
    \centering
    \caption{Fourier coefficients of Eqs.~\eqref{eq:SL_polar_Fourier}.}
    \label{tab:Fourier}
    \begin{tabular}{ll}
    \hline
    \hline
    Radial interaction function &\\
    \hline
        $a_{0,r} = 0.45579$  & \\
        $a_{1,r} = -0.97948$ & $b_{1,r}  = -1.82354$  \\
        $a_{2,r} = 0.36110$  & $b_{2,r}  = -0.07963$ \\
        $a_{3,r} = 0.29724$  & $b_{3,r}  = 0.54854$ \\
        $a_{4,r} = 0.05846$  & $b_{4,r}  = 0.09098$ \\
        $a_{5,r} = -0.11558$ & $b_{5,r}  = -0.09251$ \\
    \hline
    \hline
    Angular interaction function &\\
    $a_{0,\varphi} = 0$ & \\
    $a_{1,\varphi} = -0.00610$ & $b_{1,\varphi} = 0.31622$\\
    $a_{2,\varphi} = -0.35811$ & $b_{2,\varphi} = 0.29020$\\
    $a_{3,\varphi} = -0.25341$ & $b_{3,\varphi} = -0.0558$5\\
    $a_{4,\varphi} = -0.13541$ & $b_{4,\varphi} = 0.00799$\\
    $a_{5,\varphi} = -0.07183$ & $b_{5,\varphi} = 0.00425$\\
    \hline
    \end{tabular}
\end{table}

\section{\label{app:reduction}Model reduction to phase difference}

In order to simplify the numerical analysis, we rewrite the model~(\ref{eq:SL_polar_general}) with $N=4$ in terms of phase differences relative to oscillator~$4$: $\hat\varphi_j = \varphi_j - \varphi_4$. We approximate the terms in the interaction functions of Eq.~(\ref{eq:SL_polar_general}) as $\varphi_{j+1}(t-\tau) - \varphi_j(t) = \hat\varphi_{j+1}(t-\tau) + \varphi_4(t-\tau) - \hat\varphi_j(t) - \varphi_4(t) \approx \hat\varphi_{j+1}(t-\tau) - \hat\varphi_j(t) - \dot{\varphi}_4(t)\tau$, and rewrite the model as

\begin{subequations}
\label{eq:SL_polar_diff}
\begin{align}
  \dot{r}_j(t) =& \left[\lambda - r_j(t)^2\right]r_j(t) \\\nonumber
                &+ K r_{j+1}(t-\tau)\, H_r\left[\hat\varphi_{j+1}(t-\tau)-\hat\varphi_j(t) - \dot{\varphi}_4(t)\tau\right],\\ 
                &j=1,\dots,4\\
  \dot{\hat\varphi}_j(t) =& \omega - \gamma r_j(t)^2 \\\nonumber
                &+ K \frac{r_{j+1}(t-\tau)}{r_j(t)}\, H_\varphi\left[\hat\varphi_{j+1}(t-\tau)-\hat\varphi_j(t) - \dot{\varphi}_4(t)\tau\right]\\\nonumber
                & - \dot{\varphi}_4(t), ~j=1,2,3.
\end{align}
\end{subequations}

Equation~(\ref{eq:SL_polar_diff}b) now represents the phases of oscillators~$1$, $2$ and $3$ within the rotating frame of oscillator~$4$. Therefore, the system only requires $7$ dimensions (i.e. $j=1,\dots,4$ for Eq.~(\ref{eq:SL_polar_diff}a) and $j=1,2,3$ for Eq.~(\ref{eq:SL_polar_diff}b)). Note that when $j=3$, the term $\hat\varphi_{j+1}(t-\tau)$ in Eq.~(\ref{eq:SL_polar_diff}b) becomes zero. 

Of course, there is also the variable $\dot{\varphi}_4(t)$, which must be solved for every time Eqs.~(\ref{eq:SL_polar_diff}) are evaluated. This is done by solving the nonlinear implicit equation
\begin{subequations}
\label{eq:phi4}
\begin{align}
  0 =& \omega - \gamma r_4(t)^2 \\\nonumber
                &+ K \frac{r_1(t-\tau)}{r_4(t)}\, H_\varphi\left[\hat\varphi_{1}(t-\tau)- \dot{\varphi}_4(t)\tau\right]
                - \dot{\varphi}_4(t).
\end{align}
\end{subequations}
Since \textsc{DDE-Biftool} runs in \textsc{Matlab}, we use the function \textit{fzero} to solve Eq.~(\ref{eq:phi4}).

This reduction raises the question of whether the above approximation $\varphi_4(t-\tau) - \varphi_4(t) \approx - \dot{\varphi}_4(t)\tau$ has a significant effect on the stability of the solutions. Therefore, we calculated a sample of the above results without the reduction, and confirm no noticeable difference between the two sets of results. As an additional check, we also confirm the stability of various solutions by means of numerical simulation.

% \nocite{*}
% \bibliography{aipsamp}% Produces the bibliography via BibTeX.
\bibliography{ref.bib}% Produces the bibliography via BibTeX.

\begin{thebibliography}{35}%
\makeatletter
\providecommand \@ifxundefined [1]{%
 \@ifx{#1\undefined}
}%
\providecommand \@ifnum [1]{%
 \ifnum #1\expandafter \@firstoftwo
 \else \expandafter \@secondoftwo
 \fi
}%
\providecommand \@ifx [1]{%
 \ifx #1\expandafter \@firstoftwo
 \else \expandafter \@secondoftwo
 \fi
}%
\providecommand \natexlab [1]{#1}%
\providecommand \enquote  [1]{``#1''}%
\providecommand \bibnamefont  [1]{#1}%
\providecommand \bibfnamefont [1]{#1}%
\providecommand \citenamefont [1]{#1}%
\providecommand \href@noop [0]{\@secondoftwo}%
\providecommand \href [0]{\begingroup \@sanitize@url \@href}%
\providecommand \@href[1]{\@@startlink{#1}\@@href}%
\providecommand \@@href[1]{\endgroup#1\@@endlink}%
\providecommand \@sanitize@url [0]{\catcode `\\12\catcode `\$12\catcode
  `\&12\catcode `\#12\catcode `\^12\catcode `\_12\catcode `\%12\relax}%
\providecommand \@@startlink[1]{}%
\providecommand \@@endlink[0]{}%
\providecommand \url  [0]{\begingroup\@sanitize@url \@url }%
\providecommand \@url [1]{\endgroup\@href {#1}{\urlprefix }}%
\providecommand \urlprefix  [0]{URL }%
\providecommand \Eprint [0]{\href }%
\providecommand \doibase [0]{http://dx.doi.org/}%
\providecommand \selectlanguage [0]{\@gobble}%
\providecommand \bibinfo  [0]{\@secondoftwo}%
\providecommand \bibfield  [0]{\@secondoftwo}%
\providecommand \translation [1]{[#1]}%
\providecommand \BibitemOpen [0]{}%
\providecommand \bibitemStop [0]{}%
\providecommand \bibitemNoStop [0]{.\EOS\space}%
\providecommand \EOS [0]{\spacefactor3000\relax}%
\providecommand \BibitemShut  [1]{\csname bibitem#1\endcsname}%
\let\auto@bib@innerbib\@empty
%</preamble>
\bibitem [{\citenamefont {Juang}\ and\ \citenamefont
  {Liang}(2014)}]{juang2014}%
  \BibitemOpen
  \bibfield  {author} {\bibinfo {author} {\bibfnamefont {J.}~\bibnamefont
  {Juang}}\ and\ \bibinfo {author} {\bibfnamefont {Y.-H.}\ \bibnamefont
  {Liang}},\ }\href@noop {} {\bibfield  {journal} {\bibinfo  {journal} {Chaos}\
  }\textbf {\bibinfo {volume} {24}},\ \bibinfo {pages} {013110} (\bibinfo
  {year} {2014})}\BibitemShut {NoStop}%
\bibitem [{\citenamefont {Lodi}\ \emph {et~al.}(2020)\citenamefont {Lodi},
  \citenamefont {Della~Rossa}, \citenamefont {Sorrentino},\ and\ \citenamefont
  {Storace}}]{lodi2020}%
  \BibitemOpen
  \bibfield  {author} {\bibinfo {author} {\bibfnamefont {M.}~\bibnamefont
  {Lodi}}, \bibinfo {author} {\bibfnamefont {F.}~\bibnamefont {Della~Rossa}},
  \bibinfo {author} {\bibfnamefont {F.}~\bibnamefont {Sorrentino}}, \ and\
  \bibinfo {author} {\bibfnamefont {M.}~\bibnamefont {Storace}},\ }\href@noop
  {} {\bibfield  {journal} {\bibinfo  {journal} {Sci. Rep.}\ }\textbf {\bibinfo
  {volume} {10}},\ \bibinfo {pages} {16336} (\bibinfo {year}
  {2020})}\BibitemShut {NoStop}%
\bibitem [{\citenamefont {Protachevicz}\ \emph {et~al.}(2021)\citenamefont
  {Protachevicz}, \citenamefont {Hansen}, \citenamefont {Iarosz}, \citenamefont
  {Caldas}, \citenamefont {Batista},\ and\ \citenamefont
  {Kurths}}]{protachevicz2021}%
  \BibitemOpen
  \bibfield  {author} {\bibinfo {author} {\bibfnamefont {P.~R.}\ \bibnamefont
  {Protachevicz}}, \bibinfo {author} {\bibfnamefont {M.}~\bibnamefont
  {Hansen}}, \bibinfo {author} {\bibfnamefont {K.~C.}\ \bibnamefont {Iarosz}},
  \bibinfo {author} {\bibfnamefont {I.~L.}\ \bibnamefont {Caldas}}, \bibinfo
  {author} {\bibfnamefont {A.~M.}\ \bibnamefont {Batista}}, \ and\ \bibinfo
  {author} {\bibfnamefont {J.}~\bibnamefont {Kurths}},\ }\href@noop {}
  {\bibfield  {journal} {\bibinfo  {journal} {Front. Hum. Neurosci.}\ }\textbf
  {\bibinfo {volume} {15}},\ \bibinfo {pages} {663408} (\bibinfo {year}
  {2021})}\BibitemShut {NoStop}%
\bibitem [{\citenamefont {Soriano}\ \emph {et~al.}(2013)\citenamefont
  {Soriano}, \citenamefont {Garc{\'\i}a-Ojalvo}, \citenamefont {Mirasso},\ and\
  \citenamefont {Fischer}}]{soriano2013}%
  \BibitemOpen
  \bibfield  {author} {\bibinfo {author} {\bibfnamefont {M.~C.}\ \bibnamefont
  {Soriano}}, \bibinfo {author} {\bibfnamefont {J.}~\bibnamefont
  {Garc{\'\i}a-Ojalvo}}, \bibinfo {author} {\bibfnamefont {C.~R.}\ \bibnamefont
  {Mirasso}}, \ and\ \bibinfo {author} {\bibfnamefont {I.}~\bibnamefont
  {Fischer}},\ }\href@noop {} {\bibfield  {journal} {\bibinfo  {journal} {Rev.
  Mod. Phys.}\ }\textbf {\bibinfo {volume} {85}},\ \bibinfo {pages} {421}
  (\bibinfo {year} {2013})}\BibitemShut {NoStop}%
\bibitem [{\citenamefont {Han}\ \emph {et~al.}(2019)\citenamefont {Han},
  \citenamefont {Xiang},\ and\ \citenamefont {Zhang}}]{han2019}%
  \BibitemOpen
  \bibfield  {author} {\bibinfo {author} {\bibfnamefont {Y.}~\bibnamefont
  {Han}}, \bibinfo {author} {\bibfnamefont {S.}~\bibnamefont {Xiang}}, \ and\
  \bibinfo {author} {\bibfnamefont {L.}~\bibnamefont {Zhang}},\ }\href@noop {}
  {\bibfield  {journal} {\bibinfo  {journal} {Opt. Commun.}\ }\textbf {\bibinfo
  {volume} {445}},\ \bibinfo {pages} {262} (\bibinfo {year}
  {2019})}\BibitemShut {NoStop}%
\bibitem [{\citenamefont {Schnitzler}\ and\ \citenamefont
  {Gross}(2005)}]{schnitzler2005}%
  \BibitemOpen
  \bibfield  {author} {\bibinfo {author} {\bibfnamefont {A.}~\bibnamefont
  {Schnitzler}}\ and\ \bibinfo {author} {\bibfnamefont {J.}~\bibnamefont
  {Gross}},\ }\href@noop {} {\bibfield  {journal} {\bibinfo  {journal} {Nat.
  Rev. Neurosci.}\ }\textbf {\bibinfo {volume} {6}},\ \bibinfo {pages} {285}
  (\bibinfo {year} {2005})}\BibitemShut {NoStop}%
\bibitem [{\citenamefont {Motter}\ \emph {et~al.}(2013)\citenamefont {Motter},
  \citenamefont {Myers}, \citenamefont {Anghel},\ and\ \citenamefont
  {Nishikawa}}]{motter2013}%
  \BibitemOpen
  \bibfield  {author} {\bibinfo {author} {\bibfnamefont {A.~E.}\ \bibnamefont
  {Motter}}, \bibinfo {author} {\bibfnamefont {S.~A.}\ \bibnamefont {Myers}},
  \bibinfo {author} {\bibfnamefont {M.}~\bibnamefont {Anghel}}, \ and\ \bibinfo
  {author} {\bibfnamefont {T.}~\bibnamefont {Nishikawa}},\ }\href@noop {}
  {\bibfield  {journal} {\bibinfo  {journal} {Nat. Phys.}\ }\textbf {\bibinfo
  {volume} {9}},\ \bibinfo {pages} {191} (\bibinfo {year} {2013})}\BibitemShut
  {NoStop}%
\bibitem [{\citenamefont {Pecora}\ \emph {et~al.}(2014)\citenamefont {Pecora},
  \citenamefont {Sorrentino}, \citenamefont {Hagerstrom}, \citenamefont
  {Murphy},\ and\ \citenamefont {Roy}}]{pecora14}%
  \BibitemOpen
  \bibfield  {author} {\bibinfo {author} {\bibfnamefont {L.~M.}\ \bibnamefont
  {Pecora}}, \bibinfo {author} {\bibfnamefont {F.}~\bibnamefont {Sorrentino}},
  \bibinfo {author} {\bibfnamefont {A.~M.}\ \bibnamefont {Hagerstrom}},
  \bibinfo {author} {\bibfnamefont {T.~E.}\ \bibnamefont {Murphy}}, \ and\
  \bibinfo {author} {\bibfnamefont {R.}~\bibnamefont {Roy}},\ }\href@noop {}
  {\bibfield  {journal} {\bibinfo  {journal} {Nat. Commun.}\ }\textbf {\bibinfo
  {volume} {5}},\ \bibinfo {pages} {4079} (\bibinfo {year} {2014})}\BibitemShut
  {NoStop}%
\bibitem [{\citenamefont {MacArthur}\ \emph {et~al.}(2008)\citenamefont
  {MacArthur}, \citenamefont {S{\'a}nchez-Garc{\'\i}a},\ and\ \citenamefont
  {Anderson}}]{macarthur08}%
  \BibitemOpen
  \bibfield  {author} {\bibinfo {author} {\bibfnamefont {B.~D.}\ \bibnamefont
  {MacArthur}}, \bibinfo {author} {\bibfnamefont {R.~J.}\ \bibnamefont
  {S{\'a}nchez-Garc{\'\i}a}}, \ and\ \bibinfo {author} {\bibfnamefont {J.~W.}\
  \bibnamefont {Anderson}},\ }\href@noop {} {\bibfield  {journal} {\bibinfo
  {journal} {Discrete Applied Mathematics}\ }\textbf {\bibinfo {volume}
  {156}},\ \bibinfo {pages} {3525} (\bibinfo {year} {2008})}\BibitemShut
  {NoStop}%
\bibitem [{\citenamefont {Skardal}(2019)}]{skardal19}%
  \BibitemOpen
  \bibfield  {author} {\bibinfo {author} {\bibfnamefont {P.~S.}\ \bibnamefont
  {Skardal}},\ }\href@noop {} {\bibfield  {journal} {\bibinfo  {journal} {The
  European Physical Journal B}\ }\textbf {\bibinfo {volume} {92}},\ \bibinfo
  {pages} {1} (\bibinfo {year} {2019})}\BibitemShut {NoStop}%
\bibitem [{\citenamefont {Chossat}\ and\ \citenamefont
  {Lauterbach}(2000)}]{chossat00}%
  \BibitemOpen
  \bibfield  {author} {\bibinfo {author} {\bibfnamefont {P.}~\bibnamefont
  {Chossat}}\ and\ \bibinfo {author} {\bibfnamefont {R.}~\bibnamefont
  {Lauterbach}},\ }\href@noop {} {\emph {\bibinfo {title} {Methods in
  equivariant bifurcations and dynamical systems}}},\ Vol.~\bibinfo {volume}
  {15}\ (\bibinfo  {publisher} {World Scientific Publishing Company},\ \bibinfo
  {year} {2000})\BibitemShut {NoStop}%
\bibitem [{\citenamefont {Golubitsky}\ \emph {et~al.}(2012)\citenamefont
  {Golubitsky}, \citenamefont {Stewart},\ and\ \citenamefont
  {Schaeffer}}]{golubitsky12}%
  \BibitemOpen
  \bibfield  {author} {\bibinfo {author} {\bibfnamefont {M.}~\bibnamefont
  {Golubitsky}}, \bibinfo {author} {\bibfnamefont {I.}~\bibnamefont {Stewart}},
  \ and\ \bibinfo {author} {\bibfnamefont {D.~G.}\ \bibnamefont {Schaeffer}},\
  }\href@noop {} {\emph {\bibinfo {title} {Singularities and Groups in
  Bifurcation Theory: Volume II}}},\ Vol.~\bibinfo {volume} {69}\ (\bibinfo
  {publisher} {Springer Science \& Business Media},\ \bibinfo {year}
  {2012})\BibitemShut {NoStop}%
\bibitem [{\citenamefont {Pietras}\ and\ \citenamefont
  {Daffertshofer}(2019)}]{pietras19}%
  \BibitemOpen
  \bibfield  {author} {\bibinfo {author} {\bibfnamefont {B.}~\bibnamefont
  {Pietras}}\ and\ \bibinfo {author} {\bibfnamefont {A.}~\bibnamefont
  {Daffertshofer}},\ }\href@noop {} {\bibfield  {journal} {\bibinfo  {journal}
  {Physics Reports}\ }\textbf {\bibinfo {volume} {819}},\ \bibinfo {pages} {1}
  (\bibinfo {year} {2019})}\BibitemShut {NoStop}%
\bibitem [{\citenamefont {Nicosia}\ \emph {et~al.}(2013)\citenamefont
  {Nicosia}, \citenamefont {Valencia}, \citenamefont {Chavez}, \citenamefont
  {D{\'\i}az-Guilera},\ and\ \citenamefont {Latora}}]{nicosia13}%
  \BibitemOpen
  \bibfield  {author} {\bibinfo {author} {\bibfnamefont {V.}~\bibnamefont
  {Nicosia}}, \bibinfo {author} {\bibfnamefont {M.}~\bibnamefont {Valencia}},
  \bibinfo {author} {\bibfnamefont {M.}~\bibnamefont {Chavez}}, \bibinfo
  {author} {\bibfnamefont {A.}~\bibnamefont {D{\'\i}az-Guilera}}, \ and\
  \bibinfo {author} {\bibfnamefont {V.}~\bibnamefont {Latora}},\ }\href@noop {}
  {\bibfield  {journal} {\bibinfo  {journal} {Physical review letters}\
  }\textbf {\bibinfo {volume} {110}},\ \bibinfo {pages} {174102} (\bibinfo
  {year} {2013})}\BibitemShut {NoStop}%
\bibitem [{\citenamefont {Schneider}(2013)}]{schneider13}%
  \BibitemOpen
  \bibfield  {author} {\bibinfo {author} {\bibfnamefont {I.}~\bibnamefont
  {Schneider}},\ }\href@noop {} {\bibfield  {journal} {\bibinfo  {journal}
  {Philosophical Transactions of the Royal Society A: Mathematical, Physical
  and Engineering Sciences}\ }\textbf {\bibinfo {volume} {371}},\ \bibinfo
  {pages} {20120472} (\bibinfo {year} {2013})}\BibitemShut {NoStop}%
\bibitem [{\citenamefont {Schneider}\ and\ \citenamefont
  {Bosewitz}(2016)}]{schneider16}%
  \BibitemOpen
  \bibfield  {author} {\bibinfo {author} {\bibfnamefont {I.}~\bibnamefont
  {Schneider}}\ and\ \bibinfo {author} {\bibfnamefont {M.}~\bibnamefont
  {Bosewitz}},\ }\href@noop {} {\bibfield  {journal} {\bibinfo  {journal}
  {Disc. Cont. Dyn. Syst. A}\ }\textbf {\bibinfo {volume} {36}},\ \bibinfo
  {pages} {451} (\bibinfo {year} {2016})}\BibitemShut {NoStop}%
\bibitem [{\citenamefont {Collins}\ and\ \citenamefont
  {Stewart}(1993{\natexlab{a}})}]{collins93a}%
  \BibitemOpen
  \bibfield  {author} {\bibinfo {author} {\bibfnamefont {J.}~\bibnamefont
  {Collins}}\ and\ \bibinfo {author} {\bibfnamefont {I.}~\bibnamefont
  {Stewart}},\ }\href@noop {} {\bibfield  {journal} {\bibinfo  {journal}
  {Biological cybernetics}\ }\textbf {\bibinfo {volume} {68}},\ \bibinfo
  {pages} {287} (\bibinfo {year} {1993}{\natexlab{a}})}\BibitemShut {NoStop}%
\bibitem [{\citenamefont {Collins}\ and\ \citenamefont
  {Stewart}(1993{\natexlab{b}})}]{collins93b}%
  \BibitemOpen
  \bibfield  {author} {\bibinfo {author} {\bibfnamefont {J.~J.}\ \bibnamefont
  {Collins}}\ and\ \bibinfo {author} {\bibfnamefont {I.~N.}\ \bibnamefont
  {Stewart}},\ }\href@noop {} {\bibfield  {journal} {\bibinfo  {journal}
  {Journal of Nonlinear science}\ }\textbf {\bibinfo {volume} {3}},\ \bibinfo
  {pages} {349} (\bibinfo {year} {1993}{\natexlab{b}})}\BibitemShut {NoStop}%
\bibitem [{\citenamefont {Atay}(2010)}]{atay2010}%
  \BibitemOpen
  \bibfield  {author} {\bibinfo {author} {\bibfnamefont {F.~M.}\ \bibnamefont
  {Atay}},\ }\href@noop {} {\emph {\bibinfo {title} {Complex time-delay
  systems: theory and applications}}}\ (\bibinfo  {publisher} {Springer},\
  \bibinfo {year} {2010})\BibitemShut {NoStop}%
\bibitem [{\citenamefont {Zakharova}\ \emph {et~al.}(2013)\citenamefont
  {Zakharova}, \citenamefont {Schneider}, \citenamefont {Kyrychko},
  \citenamefont {Blyuss}, \citenamefont {Koseska}, \citenamefont {Fiedler},\
  and\ \citenamefont {Sch{\"o}ll}}]{zakharova13}%
  \BibitemOpen
  \bibfield  {author} {\bibinfo {author} {\bibfnamefont {A.}~\bibnamefont
  {Zakharova}}, \bibinfo {author} {\bibfnamefont {I.}~\bibnamefont
  {Schneider}}, \bibinfo {author} {\bibfnamefont {Y.}~\bibnamefont {Kyrychko}},
  \bibinfo {author} {\bibfnamefont {K.}~\bibnamefont {Blyuss}}, \bibinfo
  {author} {\bibfnamefont {A.}~\bibnamefont {Koseska}}, \bibinfo {author}
  {\bibfnamefont {B.}~\bibnamefont {Fiedler}}, \ and\ \bibinfo {author}
  {\bibfnamefont {E.}~\bibnamefont {Sch{\"o}ll}},\ }\href@noop {} {\bibfield
  {journal} {\bibinfo  {journal} {Europhysics Letters}\ }\textbf {\bibinfo
  {volume} {104}},\ \bibinfo {pages} {50004} (\bibinfo {year}
  {2013})}\BibitemShut {NoStop}%
\bibitem [{\citenamefont {Erneux}\ \emph {et~al.}(2017)\citenamefont {Erneux},
  \citenamefont {Javaloyes}, \citenamefont {Wolfrum},\ and\ \citenamefont
  {Yanchuk}}]{erneux2017}%
  \BibitemOpen
  \bibfield  {author} {\bibinfo {author} {\bibfnamefont {T.}~\bibnamefont
  {Erneux}}, \bibinfo {author} {\bibfnamefont {J.}~\bibnamefont {Javaloyes}},
  \bibinfo {author} {\bibfnamefont {M.}~\bibnamefont {Wolfrum}}, \ and\
  \bibinfo {author} {\bibfnamefont {S.}~\bibnamefont {Yanchuk}},\ }\href@noop
  {} {\enquote {\bibinfo {title} {Introduction to focus issue: Time-delay
  dynamics},}\ } (\bibinfo {year} {2017})\BibitemShut {NoStop}%
\bibitem [{\citenamefont {Otto}\ \emph {et~al.}(2019)\citenamefont {Otto},
  \citenamefont {Just},\ and\ \citenamefont {Radons}}]{otto2019}%
  \BibitemOpen
  \bibfield  {author} {\bibinfo {author} {\bibfnamefont {A.}~\bibnamefont
  {Otto}}, \bibinfo {author} {\bibfnamefont {W.}~\bibnamefont {Just}}, \ and\
  \bibinfo {author} {\bibfnamefont {G.}~\bibnamefont {Radons}},\ }\href@noop {}
  {\bibfield  {journal} {\bibinfo  {journal} {Philos. Trans. Royal Soc. A}\
  }\textbf {\bibinfo {volume} {377}},\ \bibinfo {pages} {20180389} (\bibinfo
  {year} {2019})}\BibitemShut {NoStop}%
\bibitem [{\citenamefont {Calleja}\ \emph {et~al.}(2017)\citenamefont
  {Calleja}, \citenamefont {Humphries},\ and\ \citenamefont
  {Krauskopf}}]{calleja2017}%
  \BibitemOpen
  \bibfield  {author} {\bibinfo {author} {\bibfnamefont {R.~C.}\ \bibnamefont
  {Calleja}}, \bibinfo {author} {\bibfnamefont {A.}~\bibnamefont {Humphries}},
  \ and\ \bibinfo {author} {\bibfnamefont {B.}~\bibnamefont {Krauskopf}},\
  }\href@noop {} {\bibfield  {journal} {\bibinfo  {journal} {SIAM J. Appl. Dyn.
  Syst.}\ }\textbf {\bibinfo {volume} {16}},\ \bibinfo {pages} {1474} (\bibinfo
  {year} {2017})}\BibitemShut {NoStop}%
\bibitem [{\citenamefont {Yanchuk}\ \emph {et~al.}(2004)\citenamefont
  {Yanchuk}, \citenamefont {Schneider},\ and\ \citenamefont
  {Recke}}]{yanchuk2004dynamics}%
  \BibitemOpen
  \bibfield  {author} {\bibinfo {author} {\bibfnamefont {S.}~\bibnamefont
  {Yanchuk}}, \bibinfo {author} {\bibfnamefont {K.~R.}\ \bibnamefont
  {Schneider}}, \ and\ \bibinfo {author} {\bibfnamefont {L.}~\bibnamefont
  {Recke}},\ }\href@noop {} {\bibfield  {journal} {\bibinfo  {journal} {Phys.
  Rev.~E}\ }\textbf {\bibinfo {volume} {69}},\ \bibinfo {pages} {056221}
  (\bibinfo {year} {2004})}\BibitemShut {NoStop}%
\bibitem [{\citenamefont {Erzgr{\"a}ber}\ \emph {et~al.}(2006)\citenamefont
  {Erzgr{\"a}ber}, \citenamefont {Krauskopf},\ and\ \citenamefont
  {Lenstra}}]{erzgraber2006compound}%
  \BibitemOpen
  \bibfield  {author} {\bibinfo {author} {\bibfnamefont {H.}~\bibnamefont
  {Erzgr{\"a}ber}}, \bibinfo {author} {\bibfnamefont {B.}~\bibnamefont
  {Krauskopf}}, \ and\ \bibinfo {author} {\bibfnamefont {D.}~\bibnamefont
  {Lenstra}},\ }\href@noop {} {\bibfield  {journal} {\bibinfo  {journal} {SIAM
  J. Appl. Dyn. Syst.}\ }\textbf {\bibinfo {volume} {5}},\ \bibinfo {pages}
  {30} (\bibinfo {year} {2006})}\BibitemShut {NoStop}%
\bibitem [{\citenamefont {Clerkin}\ \emph {et~al.}(2014)\citenamefont
  {Clerkin}, \citenamefont {O'Brien},\ and\ \citenamefont
  {Amann}}]{clerkin2014multistabilities}%
  \BibitemOpen
  \bibfield  {author} {\bibinfo {author} {\bibfnamefont {E.}~\bibnamefont
  {Clerkin}}, \bibinfo {author} {\bibfnamefont {S.}~\bibnamefont {O'Brien}}, \
  and\ \bibinfo {author} {\bibfnamefont {A.}~\bibnamefont {Amann}},\
  }\href@noop {} {\bibfield  {journal} {\bibinfo  {journal} {Phys. Rev.~E}\
  }\textbf {\bibinfo {volume} {89}},\ \bibinfo {pages} {032919} (\bibinfo
  {year} {2014})}\BibitemShut {NoStop}%
\bibitem [{\citenamefont {Della~Rossa}\ \emph {et~al.}(2020)\citenamefont
  {Della~Rossa}, \citenamefont {Pecora}, \citenamefont {Blaha}, \citenamefont
  {Shirin}, \citenamefont {Klickstein},\ and\ \citenamefont
  {Sorrentino}}]{della2020symmetries}%
  \BibitemOpen
  \bibfield  {author} {\bibinfo {author} {\bibfnamefont {F.}~\bibnamefont
  {Della~Rossa}}, \bibinfo {author} {\bibfnamefont {L.}~\bibnamefont {Pecora}},
  \bibinfo {author} {\bibfnamefont {K.}~\bibnamefont {Blaha}}, \bibinfo
  {author} {\bibfnamefont {A.}~\bibnamefont {Shirin}}, \bibinfo {author}
  {\bibfnamefont {I.}~\bibnamefont {Klickstein}}, \ and\ \bibinfo {author}
  {\bibfnamefont {F.}~\bibnamefont {Sorrentino}},\ }\href@noop {} {\bibfield
  {journal} {\bibinfo  {journal} {Nat. Commun.}\ }\textbf {\bibinfo {volume}
  {11}},\ \bibinfo {pages} {3179} (\bibinfo {year} {2020})}\BibitemShut
  {NoStop}%
\bibitem [{\citenamefont {Blaha}\ \emph {et~al.}(2016)\citenamefont {Blaha},
  \citenamefont {Burrus}, \citenamefont {Orozco-Mora}, \citenamefont
  {Ruiz-Beltr{\'a}n}, \citenamefont {Siddique}, \citenamefont {Hatamipour},\
  and\ \citenamefont {Sorrentino}}]{blaha2016symmetry}%
  \BibitemOpen
  \bibfield  {author} {\bibinfo {author} {\bibfnamefont {K.}~\bibnamefont
  {Blaha}}, \bibinfo {author} {\bibfnamefont {R.~J.}\ \bibnamefont {Burrus}},
  \bibinfo {author} {\bibfnamefont {J.~L.}\ \bibnamefont {Orozco-Mora}},
  \bibinfo {author} {\bibfnamefont {E.}~\bibnamefont {Ruiz-Beltr{\'a}n}},
  \bibinfo {author} {\bibfnamefont {A.~B.}\ \bibnamefont {Siddique}}, \bibinfo
  {author} {\bibfnamefont {V.}~\bibnamefont {Hatamipour}}, \ and\ \bibinfo
  {author} {\bibfnamefont {F.}~\bibnamefont {Sorrentino}},\ }\href@noop {}
  {\bibfield  {journal} {\bibinfo  {journal} {Chaos}\ }\textbf {\bibinfo
  {volume} {26}},\ \bibinfo {pages} {116307} (\bibinfo {year}
  {2016})}\BibitemShut {NoStop}%
\bibitem [{\citenamefont {Blaha}\ \emph {et~al.}(2013)\citenamefont {Blaha},
  \citenamefont {Lehnert}, \citenamefont {Keane}, \citenamefont {Dahms},
  \citenamefont {H{\"o}vel}, \citenamefont {Sch{\"o}ll},\ and\ \citenamefont
  {Hudson}}]{BLA13}%
  \BibitemOpen
  \bibfield  {author} {\bibinfo {author} {\bibfnamefont {K.}~\bibnamefont
  {Blaha}}, \bibinfo {author} {\bibfnamefont {J.}~\bibnamefont {Lehnert}},
  \bibinfo {author} {\bibfnamefont {A.}~\bibnamefont {Keane}}, \bibinfo
  {author} {\bibfnamefont {T.}~\bibnamefont {Dahms}}, \bibinfo {author}
  {\bibfnamefont {P.}~\bibnamefont {H{\"o}vel}}, \bibinfo {author}
  {\bibfnamefont {E.}~\bibnamefont {Sch{\"o}ll}}, \ and\ \bibinfo {author}
  {\bibfnamefont {J.~L.}\ \bibnamefont {Hudson}},\ }\href@noop {} {\bibfield
  {journal} {\bibinfo  {journal} {Phys. Rev.~E}\ }\textbf {\bibinfo {volume}
  {88}},\ \bibinfo {pages} {062915} (\bibinfo {year} {2013})}\BibitemShut
  {NoStop}%
\bibitem [{\citenamefont {Rusin}\ \emph {et~al.}(2010)\citenamefont {Rusin},
  \citenamefont {Kori}, \citenamefont {Kiss},\ and\ \citenamefont
  {Hudson}}]{rusin2010synchronization}%
  \BibitemOpen
  \bibfield  {author} {\bibinfo {author} {\bibfnamefont {C.~G.}\ \bibnamefont
  {Rusin}}, \bibinfo {author} {\bibfnamefont {H.}~\bibnamefont {Kori}},
  \bibinfo {author} {\bibfnamefont {I.~Z.}\ \bibnamefont {Kiss}}, \ and\
  \bibinfo {author} {\bibfnamefont {J.~L.}\ \bibnamefont {Hudson}},\
  }\href@noop {} {\bibfield  {journal} {\bibinfo  {journal} {Philos. Trans.
  Royal Soc. A}\ }\textbf {\bibinfo {volume} {368}},\ \bibinfo {pages} {2189}
  (\bibinfo {year} {2010})}\BibitemShut {NoStop}%
\bibitem [{\citenamefont {Balanov}\ \emph {et~al.}(2006)\citenamefont
  {Balanov}, \citenamefont {Krawcewicz},\ and\ \citenamefont
  {Steinlein}}]{balanov2006applied}%
  \BibitemOpen
  \bibfield  {author} {\bibinfo {author} {\bibfnamefont {Z.}~\bibnamefont
  {Balanov}}, \bibinfo {author} {\bibfnamefont {W.}~\bibnamefont {Krawcewicz}},
  \ and\ \bibinfo {author} {\bibfnamefont {H.}~\bibnamefont {Steinlein}},\
  }\href@noop {} {\emph {\bibinfo {title} {Applied equivariant degree}}},\
  Vol.~\bibinfo {volume} {1}\ (\bibinfo  {publisher} {American Institute of
  Mathematical Sciences Springfield},\ \bibinfo {year} {2006})\BibitemShut
  {NoStop}%
\bibitem [{\citenamefont {Engelborghs}\ \emph {et~al.}(2000)\citenamefont
  {Engelborghs}, \citenamefont {Luzyanina},\ and\ \citenamefont
  {Samaey}}]{engelborghs00}%
  \BibitemOpen
  \bibfield  {author} {\bibinfo {author} {\bibfnamefont {K.}~\bibnamefont
  {Engelborghs}}, \bibinfo {author} {\bibfnamefont {T.}~\bibnamefont
  {Luzyanina}}, \ and\ \bibinfo {author} {\bibfnamefont {G.}~\bibnamefont
  {Samaey}},\ }\href@noop {} {\bibfield  {journal} {\bibinfo  {journal} {TW
  Report}\ }\textbf {\bibinfo {volume} {305}} (\bibinfo {year}
  {2000})}\BibitemShut {NoStop}%
\bibitem [{\citenamefont {Sieber}\ \emph {et~al.}(2014)\citenamefont {Sieber},
  \citenamefont {Engelborghs}, \citenamefont {Luzyanina}, \citenamefont
  {Samaey},\ and\ \citenamefont {Roose}}]{SIE14a}%
  \BibitemOpen
  \bibfield  {author} {\bibinfo {author} {\bibfnamefont {J.}~\bibnamefont
  {Sieber}}, \bibinfo {author} {\bibfnamefont {K.}~\bibnamefont {Engelborghs}},
  \bibinfo {author} {\bibfnamefont {T.}~\bibnamefont {Luzyanina}}, \bibinfo
  {author} {\bibfnamefont {G.}~\bibnamefont {Samaey}}, \ and\ \bibinfo {author}
  {\bibfnamefont {D.}~\bibnamefont {Roose}},\ }\href@noop {} {\emph {\bibinfo
  {title} {{DDE-BIFTOOL} Manual - Bifurcation analysis of delay differential
  equations}}} (\bibinfo {year} {2014}),\ \Eprint
  {http://arxiv.org/abs/1406.7144} {arXiv:1406.7144 [math.DS]} \BibitemShut
  {NoStop}%
\bibitem [{\citenamefont {Pyragas}(1992)}]{pyragas92}%
  \BibitemOpen
  \bibfield  {author} {\bibinfo {author} {\bibfnamefont {K.}~\bibnamefont
  {Pyragas}},\ }\href@noop {} {\bibfield  {journal} {\bibinfo  {journal}
  {Physics letters A}\ }\textbf {\bibinfo {volume} {170}},\ \bibinfo {pages}
  {421} (\bibinfo {year} {1992})}\BibitemShut {NoStop}%
\bibitem [{\citenamefont {Zhai}\ \emph {et~al.}(2008)\citenamefont {Zhai},
  \citenamefont {Kiss},\ and\ \citenamefont {Hudson}}]{ZHA08}%
  \BibitemOpen
  \bibfield  {author} {\bibinfo {author} {\bibfnamefont {Y.}~\bibnamefont
  {Zhai}}, \bibinfo {author} {\bibfnamefont {I.~Z.}\ \bibnamefont {Kiss}}, \
  and\ \bibinfo {author} {\bibfnamefont {J.~L.}\ \bibnamefont {Hudson}},\
  }\href@noop {} {\bibfield  {journal} {\bibinfo  {journal} {Ind. Eng. Chem.
  Res.}\ }\textbf {\bibinfo {volume} {47}},\ \bibinfo {pages} {3502} (\bibinfo
  {year} {2008})}\BibitemShut {NoStop}%
\end{thebibliography}%

\end{document}